\documentclass[]{elsarticle}
\usepackage{graphicx}
\usepackage{color}
\usepackage{tikz-cd}
\usepackage{amsmath}
\usepackage{amssymb}
\usepackage{amsthm}
\usepackage{stmaryrd}
\usepackage{wasysym}
\usepackage{booktabs}
\usepackage{units}
\usepackage{mathrsfs}
\usepackage{marvosym}

\theoremstyle{theorem}
  \newtheorem{thm}{Theorem}[section]
  \newtheorem{cor}[thm]{Corollary}
  \newtheorem{lem}[thm]{Lemma}
  \newtheorem{propn}[thm]{Proposition}
  
\theoremstyle{definition}
  \newtheorem{defn}[thm]{Definition}
  \newtheorem{rem}[thm]{Remark}
  
  \newtheorem{exm}[thm]{Example}

  \DeclareFontFamily{U}{mathc}{}
  \DeclareFontShape{U}{mathc}{m}{it}%
  {<->s*[1.03] mathc10}{}
  \DeclareMathAlphabet{\mathfun}{U}{mathc}{m}{it}
  
\newcommand{\cat}[1]{\mathbf{#1}}

\newcommand{\fun}[1]{\mathfun{#1}}

\newcommand{\mf}[1]{\mathfrak{#1}}
\newcommand{\ms}[1]{\mathscr{#1}}
\newcommand{\mb}[1]{\mathbb{#1}}
\newcommand{\lan}[1]{\mathsf{#1}}

\newcommand{\mo}[1]{\mathfrak{#1}}
\newcommand{\model}[1]{\mo{#1}}
\newcommand{\dual}[1]{{#1}^{\partial}}

\newcommand{\ov}[1]{\overline{#1}}
\newcommand{\undto}{\mathrel{\mkern1mu\underline{\mkern-1mu \to\mkern-2mu}\mkern2mu }}
\newcommand{\undbito}{\mathrel{\mkern3.4mu\underline{\mkern-3.4mu \bito\mkern-3.4mu}\mkern3.4mu }}

\DeclareMathOperator{\Prop}{Prop}

\DeclareMathOperator{\Tbiint}{\lan{Tense}}
\newcommand{\llb}{\llbracket}
\newcommand{\rrb}{\rrbracket}
\renewcommand{\iff}{\quad\text{iff}\quad}
\newcommand{\bisim}{\rightleftharpoons}
\newcommand{\logeq}{\leftrightsquigarrow}
\newcommand{\preq}{\preccurlyeq}

\mathchardef\hyphen="2D
\newcommand{\To}{\Rightarrow}
\newcommand{\curvearrowleftup}{{\mathbin{\rotatebox[origin=c]{180}{$\curvearrowleft$}}}}
\newcommand{\gen}{\rotatebox[origin=c]{180}{\reflectbox{$\neg$}}}

\newcommand{\Nat}{\mathbb{N}}
\newcommand{\Lint}{\mathsf{Int}}
\newcommand{\Ldint}{\mathsf{Int}^{\partial}}
\newcommand{\Lbiint}{\mathsf{Bi\hyphen int}}
\newcommand{\bilogeq}{\logeq_{\Lbiint}}
\newcommand{\ilogeq}{\logeq_{\Lint}}
\newcommand{\dlogeq}{\logeq_{\Ldint}}
\newcommand{\lsem}{\llb}
\newcommand{\rsem}{\rrb}

\newcommand{\bito}{\operatorname{\,
\begin{tikzpicture}
\draw (0,0) --(1.5ex,0) -- (2.2ex,0.5ex);
\draw (1.5ex,0) -- (2.2ex,-0.5ex);
\end{tikzpicture}\,
}}

\renewcommand{\Diamond}{\mathchoice
    {\rotatebox[origin=c]{45}{\scalebox{.85}{$\Box$}}}
    {\rotatebox[origin=c]{45}{\scalebox{.85}{$\boxempty$}}}
    {\rotatebox[origin=c]{45}{\scalebox{.5}{$\Box$}}}
    {\rotatebox[origin=c]{45}{\scalebox{.4}{$\Box$}}}
  }
\newcommand{\cbox}{{\boxbar}}
\newcommand{\cdiamond}{{\mathbin{\rotatebox[origin=c]{45}{$\boxslash$}}}}

\newcommand{\tbox}{\mathchoice
    {\scalebox{.8}{$\blacksquare$}}
    {\scalebox{.8}{$\blacksquare$}}
    {\scalebox{.55}{$\blacksquare$}}
    {\scalebox{.45}{$\blacksquare$}}
  }
\newcommand{\tdiamond}{\mathchoice
    {\rotatebox[origin=c]{45}{\scalebox{.75}{$\blacksquare$}}}
    {\rotatebox[origin=c]{45}{\scalebox{.75}{$\blacksquare$}}}
    {\rotatebox[origin=c]{45}{\scalebox{.5}{$\blacksquare$}}}
    {\rotatebox[origin=c]{45}{\scalebox{.4}{$\blacksquare$}}}
  }

\begin{document}

\begin{frontmatter}

\title{Hennessy-Milner Properties via Topological Compactness}

\author{Jim de Groot (\Letter)}
\ead{jim.degroot@anu.edu.au}
\author{Dirk Pattinson}
\ead{dirk.pattinson@anu.edu.au}
\address{%
  School of Computing\\
  The Australian National University\\
  Canberra, Australia} 

\begin{abstract}
  We give Hennessy-Milner classes for
  intuitionistic, dual-intuitionistic and bi-intuitionistic logic
  interpreted in intuitionistic Kripke models,
  and generalise these results to modal (dual- and bi-)intuitionistic logics.
  Our main technical tools are image-compact and
  pre-image-compact relations that provide a semantical description
  of modal saturation properties.
\end{abstract}

\begin{keyword}
  Bisimulation
  \sep Hennessy-Milner property
  \sep Intuitionistic logic
  \sep Bi-intuitionistic logic
  \sep Modal logic 
\end{keyword}

\end{frontmatter}

\section{Introduction}

  \noindent
  Bisimulations play a crucial role in the model theory of modal
  logic as the canonical notion of \emph{semantic} equivalence:
  bisimilar worlds necessarily  satisfy precisely the same formulae.
  If the converse is also true, the (usually finitary) logical
  language is powerful enough to describe the (typically infinitary)
  semantics: this is the so-called \emph{Hennessy-Milner property}~\cite{HenMil85}.
  
  Bisimulations were introduced in \cite{Ben76} to
  characterise normal modal logic over a classical base as the bisimulation-invariant
  fragment of first-order logic. 
  Independently, they arose in the field of computer science as an
  equivalence relation between process graphs \cite{Mil80,Par81},
  and as extensional equality in 
  non-wellfounded set theory \cite{Acz88}.

  By and large, the Hennessy-Milner property is well understood for
  normal modal logic over a classical base, where it is known to hold for all modally
  saturated models, see Section 2 of \cite{BRV01}. 
  In the realm of (dual- and bi-)intuitionistic logic and their modal extensions,
  much less is known.
  Some explorations are made in \cite{Pat97} where the Hennessy-Milner property is
  established for \emph{intuitionistic} propositional logic,
  interpreted over intuitionistic Kripke models \cite{Kri65},
  and in \cite{Dav09}, where a Hennessy-Milner property is given for
  tense intuitionistic logic where all modalities are interpreted using
  a single additional relation.
  Besides, \cite{Bou04} contains Hennessy-Milner results for
  strict-weak languages, and \cite{Pre14} discusses a Hennessy-Milner result
  for unimodal extensions of positive, intuitionistic and bi-intuitionistic logic.

  In this paper we aim to derive Hennessy-Milner properties for a large variety
  of logics using the notion of \emph{image-com\-pact\-ness}.
  A relation is image-compact if its successor sets of a single points
  are compact in a topology that includes all truth sets of formulae as clopens.
  Similar methods have previously been used in the setting of normal modal logic
  over a classical base \cite{BonKwi95} and unimodal logic over a positive,
  intuitionistic and bi-intuitionistic base \cite{Pre14}.
  Our results apply to intuitionistic, dual-intuitionistic, and bi-intuitionistic
  propositional logic, as well as their extension with normal modal
  operators.
  Moreover, we can use them to obtain new
  Hennessy-Milner type results for various logics previously
  studied, notably modal intuitionistic, and tense bi-intuitionistic logic.

  Technically, we show that logical equivalence and bisimulations
  coincide for image-compact Kripke models, and obtain a (known)
  characterisation for intuitionistic propositional logic. We then
  dualise the semantics to obtain the same result for 
  \emph{dual-intuitionistic logic}, which
  is the extension of positive logic with a binary subtraction arrow $\bito$
  residuated with respect to disjunction.
  While this may seem like a mathematical curiosity at first,
  subtraction has found multiple applications. 
  In computer science it
  can be used to describe control mechanisms such as co-routines
  \cite{Cro04} and in philosophy the subtraction arrow provides a tool to reason about 
  refutation \cite{Tra17,Res97}.

  Thereafter, we merge the results for intuitionistic propositional logic
  and its dual to obtain a characterisation of bisimulation for
  bi-intuitionistic logic (which  can be viewed as the union of
  intuitionistic and dual-intuitionistic logic) in terms of logical equivalence.
  Bi-intuitionistic logic is also known as
  subtractive logic \cite{Cro04} and Heyting-Brouwer logic \cite{Rau74b},
  and was introduced by Rauszer with Kripke semantics and a Hilbert calculus \cite{Rau80}.
  We refer to \cite{GorShi20} for an excellent overview of the logic,
  that moreover clarifies some of Rauszer's confusions.
  
  In a second step, we extend the underlying propositional languages
  with modal operators that are interpreted like 
   Bo\v{z}i\'{c} and Do\v{s}en did in
  \cite{BozDos84}, where $\Box$ and $\Diamond$ are a priori
  unrelated modalities. Our approach is similar to the propositional
  case: a Hennessy-Milner theorem for intuitionistic propositional
  logic augmented with $\Box$ gives, by duality, an analogous
  theorem for dual-intuitionistic logic with $\Diamond$, and both
  can be combined to get the same for bi-intuitionistic logic,
  extended with an arbitrary number of $\Box$ and
  $\Diamond$-operators.

  Finally, we apply our results to obtain
  new Hennessy-Milner theorems for a large variety of logics studied
  in the literature. These fall into two classes: various flavours
  of intiutionistic modal logic
  \cite{Ono77,Fis81,PloSti86,WolZak98} and various flavours of tense
  bi-intuitionistic logic 
  \cite{GorPosTiu10,SteSchRye16,SanSte17}.

\paragraph{Structure of the Paper}
  In Section \ref{sec:IKM} we recall intuitionistic Kripke frames and models as semantics for
  intuitionistic, dual-intuitionistic and bi-intuitionistic logic.
  We give the definition of general frames and use these to define the notions
  of image-compactness and pre-image-compactness.
  Subsequently, in Section \ref{sec:inst}, we show how one can relate the relations of logical
  equivalence for different languages, 
  borrowing a simple observation from the theory of institutions.
  
  Bisimulations between intuitionistic Kripke models are defined in
  Section \ref{sec:non-modal}, and the notions of (pre-)image-compactness
  are shown to give rise to Hennessy-Milner type results for (bi- and dual-)intuitionistic
  logic.
  
  In Section \ref{sec:modal} we extend our scope to modal extensions of
  the previously studied logics. We give a suitable notion of frame and model and
  define bisimulations between them.
  Again, the notions of (pre-)image-compactness give rise to Hennessy-Milner results.
  We then specialise these results to obtain Hennessy-Milner theorems for
  a number of logics studied in the literature in Section \ref{sec:app}.
  
  Finally, in Section \ref{sec:ic-vs-sat} we detail how in some cases image-compactness
  coincides with notions of saturation, and in Section \ref{sec:conc} we
  suggest several avenues for further research.

\paragraph{Related Work}

  As mentioned above, in \cite{Bou04} the author proves Hennessy-Milner type theorems
  for strict-weak languages. Amongst such languages are intuitionistic logic,
  where implication is viewed as a strict arrow, dual-intuitionistic logic,
  modelling subtraction as a weak arrow, and bi-intuitionistic logic.
  In fact, the framework in {\it op.~\!cit.}~allows one to add as many such
  arrows as desired.
  The strict and weak arrows are interpreted using a relation
  in the same way implication and subtraction are interpreted
  (see Section \ref{sec:IKM} below).
  Moreover, every arrow gives rise to a box- or diamond-like modality via
  $\Box\phi := \top \to \phi$ and $\Diamond\phi := \phi \leftharpoondown \bot$,
  where $\leftharpoondown$ denotes a weak arrow.
  However, boxes and diamonds are not defined separately.
  This means that, when proving that some relation satisfies the back-and-forth
  conditions of a bisimulation, one can always
  make use of the arrows interpreted via each relation in the frame.
  This simplifies the proof of Hennessy-Milner results, because each clause
  resembles the proof of \cite[Theorem 21]{Pat97} or Theorem \ref{thm:hm-int} below,
  or its dual.
  In Section \ref{sec:modal} of the current paper, dealing with normal modal extensions
  of (bi- and dual-)intuitionistic logic, we do not have this luxury.

  In \cite{Pre14} the author considers modal extensions of positive, intuitionistic
  and bi-intuitionistic logic. Moreover, the relation used to interpret
  the modalities is not required to interact with the underlying partial order at all.
  The level of generality forces to author to obtain a Hennessy-Milner theorem
  via a duality, because the potential absence of implication or subtraction
  arrow frustrates a more direct approach like in \cite[Theorem 21]{Pat97}
  or \cite[Proposition 2.54]{BRV01}.
  By cleverly extending the duality to a dual adjunction, a slightly larger
  Hennessy-Milner class is derived.
  However, the models it contains are still based on pre-Priestley spaces.
  In our setting we begin with (bi- or dual-)intuitionistic logic,
  so that we always have an arrow in our language.
  Furthermore, the relations we use to interpret additional modal operators
  are required to satisfy certain coherence conditions with respect to the
  pre-order underlying a frame. These extra constraints allow us to derive
  a stronger Hennessy-Milner result.
  
  Finally, in \cite{Dav09} the author derives a Hennessy-Milner theorem for
  tense intuitionistic logic. This is a bit farther removed from our research,
  because the underlying intuitionistic logic is interpreted in topological spaces,
  rather than the more restrictive intuitionistic Kripke frames (= Alexandrov spaces)
  used here. We discuss this setting as a potential avenue for further research
  in the conclusion.

\paragraph{Relation to Predecessor Paper}
  The current paper is an extension of preliminary
  results reported in \cite{GroPat19}.
  Conceptually, we identify the core notion of image compactness as
  the key stepping stone in establishing Hennessy-Milner type
  theorems.
  Technically, this yields stronger results:  in {\it op.~\!cit.}, we
  have established Hennessy-Miler type theorems for
  descriptive and finite models of bi-intuitionistic logic.
  Both are special cases of
  (pre-)image-compact models. Moreover, (pre-)image compact models are closed
  under disjoint unions whence closure under disjoint unions,
  reported in {\it op.~\!cit.}, is automatic, and all results follow
  from Theorem \ref{thm:hm-bi-int} below.
  Similarly, the results from 
  Section 5 of \cite{GroPat19} about 
  descriptive and finite $\Lbiint_{\Box\Diamond}$-models are
  subsumed by Theorem \ref{thm:hm-mod-bi-int}, again noting that
  image compactness subsumes both finiteness and being descriptive.
  Finally, the treatment of bisimulations for modal and epistemic intuitionistic,
  and tense bi-intuitionistic logic is new.

\section{Intuitionistic Kripke Models and Image Compactness}\label{sec:IKM}

\noindent
  We recall the Kripke semantics of intuitionistic, dual-intuitionistic and
  bi-intuitionistic propositional logic, and introduce the semantic notion at
  the heart of our results: image-compact relations.
  Throughout the paper, we write $\Prop$ for a (possibly infinite) set of
  propositional variables.

\begin{defn}
The \emph{language} $\Lbiint(\Prop)$ of  bi-intuitionistic
propositional logic
over the set $\Prop$ of
propositional variables is given by the grammar
  $$
    \phi ::= \top \mid \bot \mid p \mid \phi \wedge \phi \mid \phi \vee \phi 
    \mid \phi \to \psi \mid \phi \bito \phi .
  $$
where $\to$ is intuitionistic implication and $\bito$ its dual,
sometimes called \emph{subtraction}.

The language $\Lint(\Prop)$ of intuitionistic propositional logic is the
set of $\bito$-free bi-intuitionistic formulae, and the language
$\Ldint(\Prop)$ consists of all implication-free formulae.
\end{defn}

\noindent
  All three languages can be interpreted over intuitionistic Kripke models.
  These are simply pre-ordered sets, i.e., sets with a reflexive and transitive
  relation on them.
  If $(X, \leq)$ is a pre-order and $a \subseteq X$ then we write
  ${\uparrow}a = \{ y \in X \mid x \leq y \text{ for some } x \in a \}$ for the upwards closure of $a$, and for $x \in X$ we abbreviate ${\uparrow}x := {\uparrow}\{ x \}$.
  The set $a$ is called an \emph{upset} if ${\uparrow}a = a$, and we write
  $\fun{Up}(X, \leq)$ for the collection of upsets of $(X, \leq)$.

\begin{defn}\label{def:interpr}
  An \emph{intuitionistic Kripke frame} is a pre-ordered set $(X, \leq)$.
  An \emph{intuitionistic Kripke model} is a triple $(X, \leq, V)$
  where $(X, \leq)$ is a pre-order, and $V: \Prop \to \fun{Up}(X, \leq)$
  is an upset-valued valuation.
  
  The \emph{truth} of bi-intuitionistic formulae in an intuitionistic
  Kripke model $\mo{M} = (X, \leq, V)$ at a world $x \in X$ is
  defined inductively by
  \begin{align*}
    \mo{M}, x \Vdash \top &\quad\text{always} \\
    \mo{M}, x \Vdash \bot &\quad\text{never} \\
    \mo{M}, x \Vdash p &\iff x \in V(p) \\
    \mo{M}, x \Vdash \phi \wedge \psi &\iff x \Vdash \phi \text{ and } x \Vdash \psi \\
    \mo{M}, x \Vdash \phi \vee   \psi &\iff x \Vdash \phi \text{ or } x \Vdash \psi \\
    \mo{M}, x \Vdash \phi \to \psi &\iff \text{for all } y \geq x,
                                 \text{ if  } y \Vdash \phi
                                 \text{ then } y \Vdash \psi \\
    \mo{M}, x \Vdash \phi \bito \psi &\iff \text{there exists } y \leq x
                                 \text{ such that } y \Vdash \phi
                                 \text{ xand } y \not\Vdash \psi.
  \end{align*}
  We write $x \bilogeq x'$ to denote that two states $x \in
  X$ and $x' \in X'$ of two intuitionistic Kripke models $\mo{M} = (X, \leq,
  V)$ and $\mo{M}' = (X', \leq', V')$ are logically equivalent with respect to
  bi-intuitionistic propositional logic, i.e.,
  \[ 
    \mo{M}, x \Vdash \phi \iff \mo{M}', x' \Vdash \phi
  \]
  for all $\phi \in \Lint$. The relations $\ilogeq$ and
  $\dlogeq$ are the relations of logical equivalence with
  respect to $\Lint$ and $\Ldint$ are defined analogously.
  In an intuitionistic Kripke model $\mo{M} = (X, \leq, V)$, we write 
  $\lsem \phi \rsem^{\mo{M}} = \lbrace x \in X \mid x \Vdash \phi \rbrace$
  for the truth set of $\phi$ in $\mo{M}$.
\end{defn}
  If we define the operators  $\undto, \undbito : \fun{Up}(X,
  \leq) \times \fun{Up}(X, \leq) \to \fun{Up}(X, \leq)$ by
  \begin{align*} 
    a \undto b &= \{ x \in X \mid \text{for all } y \in X, \text{ if } x \leq y \text{ and } y \in a \text{ then } y \in b \} 
                  \\
    a \undbito b &= \{ x \in X \mid \text{there exists } y \leq x
                                    \text{ such that } y \in a
                                    \text{ and } y \notin b \}
  \end{align*}
  then evidently
  $\lsem \phi \to \psi \rsem^{\mo{M}} = \lsem \phi \rsem^{\mo{M}} \undto \lsem \psi
  \rsem^{\mo{M}}$ and $\lsem \phi \bito \psi \rsem^{\mo{M}} = \lsem \phi \rsem^{\mo{M}}
  \undbito \lsem \psi \rsem^{\mo{M}}$ for any intuitionistic Kripke
  model $\mo{M}$.

  The logics $\Lint, \Ldint, \Lbiint$
  are sometimes interpreted over posets (rather
  than pre-orders), for example in the predecessor paper of this one
  \cite{GroPat19} and in \cite{ChaZak97}. Here, we choose the more general semantics.

  The relationship between intuitionistic and dual-intuitionistic
  logic is best clarified in terms of dual models (with reversed
  order). 

\begin{defn}\label{def:ikm-dual}
  The \emph{dual} of an intuitionistic Kripke model $\model{M} = (X, \leq, V)$
  is the model $\dual{\model{M}} = (X, \geq, \dual{V})$,
  where $\dual{V}$ is defined by $\dual{V}(p) = X \setminus V(p)$.
\end{defn}

  \noindent
  The notion of dual model is well defined, as the complement $X
  \setminus a$ of an
  upset $a$ in a pre-order $(X, \leq)$ is a downset, and hence an
  upset for the dual pre-order $(X, \geq)$.
  On the level of languages, we have a translation $(\cdot)^t:
  \Lint \to \Ldint$ such that
  $\phi$ is true at a state $x$ in a model $(X, \leq, V)$ if and only if
  its translation $\phi^t$ is \emph{false} in the \emph{dual model}. We define
  this inductively via
  \begin{align*}
    \bot^t &= \top  &
    \top^t &= \bot &
    p^t &= p \\
    && (\phi \wedge \psi)^t &= \phi^t \vee \psi^t &
    (\phi \vee \psi)^t &= \phi^t \wedge \psi^t  \\
    && (\phi \to \psi)^t &= \psi^t \bito \phi^t &
    (\phi \bito \psi)^t &= \psi^t \to \phi^t
  \end{align*}
  Clearly, $(\cdot)^t$ is an involution of $\Lbiint$ which restricts to translations
  $\Lint \to \Ldint$ and $\Ldint \to \Lint$.

\begin{lem}\label{lem:sem-trans}
  Let $\mo{M} = (X, \leq, V)$ be an intuitionistic Kripke model and
  $\phi \in \Lbiint$ be a formula. Then we have
  $$
    \mo{M}, x \Vdash \phi \iff \dual{\mo{M}}, x \not\Vdash \phi^t.
  $$
\end{lem}
\begin{proof}
  This follows from a straightforward induction.
  We showcase one of the inductive steps:
  \begin{align*}
    \mo{M}, x \Vdash \phi \to \psi
      &\iff \text{for all } y \geq x
            \text{ either } \mo{M}, y \not\Vdash \phi
            \text{ or } {M}, y \Vdash \psi \\   
      &\iff \text{for all } y \geq x
            \text{ either } \dual{\mo{M}}, y \Vdash \phi^t
            \text{ or } \dual{\mo{M}}, y \not\Vdash \psi^t \\   
      &\iff \text{there is no } y \geq x
            \text{ such that } \dual{\mo{M}}, y \Vdash \psi^t
            \text{ and } \dual{\mo{M}}, y \not\Vdash \phi^t \\
      &\iff \dual{\mo{M}}, x \not\Vdash \psi^t \bito \phi^t = (\phi \to \psi)^t
  \end{align*}
  All other cases are similar.
\end{proof}

  \noindent
  We now define image-compactness, the main technical vehicle that we use to
  establish Hennessy-Milner results in this paper.
  For this, we  augment models with a collection of 
  \emph{admissible subsets}, that is, a selection of subsets of the
  carrier that includes all truth sets. This allows us to
  topologise the model using the patch topology, and use compactness
  to get a finitary handle on the successors of any given world.

\begin{defn}\label{def:general}
  A \emph{general model} is a tuple
  $\mo{M} = (X, \leq, V, A)$ such that ${(X, \leq, V)}$ is an intuitionistic Kripke
  model, $A \subseteq \fun{Up}(X, \leq)$ is a collection of
  up-closed subsets of $(X, \leq)$ that (i)
  is closed under finite union and finite intersection,
  and (ii)
  contains $\emptyset$, $X$ and $V(p)$ 
  for every $p \in \Prop$.

  We call $\mo{M}$ a \emph{general $\Lint$-model} (resp.
  $\Ldint$-model)  if $A$ is moreover closed under $\undto$
  (resp. $\undbito$), and a \emph{general $\Lbiint$-model} if $A$
  is closed under both $\undto$ and $\undbito$.
 
 The \emph{patch topology} on a general model $\mo{M} = (X, \leq, V,
 A)$ is the topology $\tau_A$ on $X$ generated 
  by the (clopen) subbase $A \cup -A$,
  where $-A = \{ X \setminus a \mid a \in A \}$.
\end{defn}

\noindent
What will be of special interest later are the \emph{compact} subsets of a
general model $\mo{M} = (X, \leq, V, A)$. Recall that a subset $U \subseteq X$
is \emph{compact} if every open cover $(O_i)_{i \in I}$ of $U$ (that is, $U
\subseteq \bigcup \lbrace O_i \mid i \in I \rbrace$ and $O_i \in \tau_A$
for all $i \in I$) has a finite
subcover (that is, there exists a finite $J \subseteq I$ such that $U
\subseteq \bigcup \lbrace O_j \mid j \in J \rbrace$). 

In particular, if $x \in X$ is a world in a model $(X, \leq, V, A)$,
then bisimulation requires us to establish a property for \emph{all}
successors in $\mo{M}$, i.e., for the set
${\uparrow}_{\leq}x = \{ y \in X \mid x \leq y \}$.
If ${\uparrow}_{\leq}x$ is compact, this can be achieved in a finitary way.
This motivates the following definition of image-compactness.

\begin{defn}
  An intuitionistic Kripke model $(X, \leq, V)$ is called \emph{(pre-)image-compact}
  for $\lan{L}$ (where $\lan{L} \in \{ \Lint, \Ldint, \Lbiint \}$)
  if there exists a set $A$ of admissibles such that $(X, \leq, A, V)$
  is a general $\lan{L}$-model and for all $x \in X$ the set ${\uparrow}_{\leq}x$
  (resp.~${\downarrow}_{\leq}x$) is compact in the patch topology $\tau_A$.
\end{defn}

\noindent
  Observe that, like saturation, (pre-)image-compactness is a property of
  \emph{models}, rather than a property of frames.
  Furthermore, note that by definition of the patch topology, proposition letters are
  interpreted as clopen sets in this topology.
  We conclude the section with the following examples.

\begin{exm}

\begin{enumerate}
  \item A Kripke model $\mo{M} = (X, \leq, V)$ is \emph{image-finite} if the set
  $\lbrace y \in X \mid x \leq y \rbrace$ is finite for every $x \in X$.
  Clearly every image-finite Kripke model is image-compact: take $A$
  to be the collection of all upward closed subsets of $W$.
  \item Image-compact is strictly more general than image-finite.
  Consider for example $X = \mb{N} \cup \lbrace \infty \rbrace$ where
  $n \leq \infty$ for all $n \in \mb{N} \cup \{ \infty \}$
  (and $\leq$ is as usual otherwise), with
  the valuation $V(p_i) = \lbrace x \in X \mid i \leq x \rbrace$,
  for $i \in \mb{N}$.
  Clearly, this is not image-finite.
  If we take $A$ to consist of all sets of the form
  $\lbrace x \in X \mid x \geq n \rbrace$ where $n$ ranges over $\Nat$,
  then this is easily seen to be image-compact.
  \item Every descriptive intuitionistic Kripke frame \cite[Section 8.4]{ChaZak97} 
  is automatically
  image-compact. This follows because descriptive frames are
  precisely Esakia spaces \cite{Esa74}, hence
  topologically compact, and upsets of single points are closed in this topology.
  \item If $\mo{M} = (X, \leq, V, A)$ is a general model, and $\dual{\mo{M}} = (X,
  \geq, \dual{V}, \dual{A})$ is its dual where $\dual{A} = \lbrace X \setminus a \mid a
  \in A \rbrace$, then $\mo{M}$ is image-compact if and only if
  $\dual{\mo{M}}$ is pre-image-compact.
\end{enumerate}
\end{exm}

\section{Relating Logical Equivalence for Different Logics}\label{sec:inst}

\noindent
  As this paper is concerned with many different logics,
  it is useful to structure the relationships between them.
  More precisely, we will often show that the relation of logical equivalence between
  two models is a bisimulation for a certain logic.
  The following simple fact, borrowed from the theory of
  institutions \cite{Goguen:1992:IAM}, allows us to transfer such results from
  one logic to another.
  
  Let us abstractly define a \emph{semantics} for a language $\lan{L}$ to be
  a class of models $\mb{M}$ such that:
  \begin{itemize}
    \item Each $\mo{M} \in \mb{M}$ has an underlying set, denoted by $\fun{U}\mo{M}$; and
    \item Each model $\mo{M} \in \mb{M}$ comes with a theory map
          $\fun{th}_{\mo{M}} : \fun{U}\mo{M} \to \fun{P}\lan{L}$
          that sends a state $x \in \fun{U}\mo{M}$ to the collection of
          $\lan{L}$-formulae true at that state.
          ($\fun{P}\lan{L}$ denotes the powerset of $\lan{L}$.)
  \end{itemize}
  
  The collection $\mb{M}$ may be regarded as a category and
  $\fun{U}$ as a functor $\mb{M} \to \cat{Set}$ from $\mb{M}$ to the category of
  sets. However, we do not need this categorical perspective for our purposes.

\begin{exm}
  One can think of $\lan{L} = \Lint$, with $\mb{M}$ the collection of
  intuitionistic Kripke models from Definition \ref{def:interpr}.
  Then for $\mo{M} = (X, \leq, V) \in \mb{M}$,
  the underlying set is given by $\fun{U}\mo{M} = X$ and the theory map is
  induced by the interpretation from Definition \ref{def:interpr} via
  $$
    \fun{th}_{\mo{M}} : X \to \fun{P}\Lint : x \mapsto \{ \phi \in \Lint \mid x \Vdash \phi \}.
  $$
  It is easy to see that, {\it mutatis mutandis}, this yields semantics for 
  $\Ldint$ and $\Lbiint$ as well. 
\end{exm}
  
  If we have sufficient coherence between two such logic,
  then logical equivalence of one implies logical equivalence of the other.
  The next lemma describes this in detail.

\begin{lem}\label{lem:inst}
  Let $\lan{L}_1$ and $\lan{L}_2$ be two languages with semantics
  $\mb{M}_1$ and $\mb{M}_2$. Denote the underlying set of a model
  $\mo{M} \in \mb{M}_i$ by $\fun{U}_i\mo{M}$,
  and the theory of $x \in \fun{U}_i\mo{M}$ by $\fun{th}_i(\mo{M})(x)$.
  Let
  \begin{itemize}
    \item $t: \lan{L}_1 \to \lan{L}_2$ is a surjective translation from
          $\lan{L}_1$ to $\lan{L}_2$; and
    \item $r: \mb{M}_2 \to \mb{M}_1$ a transformation of models such that
          $\fun{U}_1(r\mo{M}) = \fun{U}_2\mo{M}$ for all $\mo{M} \in \mb{M}_2$.
  \end{itemize}
  Moreover, suppose that
  \begin{equation}\label{eq:inst-th}
    t(\phi) \in \fun{th}_2({\mo{M}})(x) \iff \phi \in \fun{th}_1({r\mo{M}})(x)
  \end{equation}
  for all $\phi \in \lan{L}_1$ and $\mo{M} \in \mb{M}_2$ and $x \in \fun{U}_2\mo{M}$.
  Then we have
  \begin{equation}\label{eq:inst-le}
    \fun{th}_2({\mo{M}})(x) = \fun{th}_2({\mo{M}'})(y)
      \iff \fun{th}_1({r\mo{M}})(x) = \fun{th}_1({r\mo{M}'})(y)
  \end{equation}
  for all $\mo{M}, \mo{M}' \in \mb{M}_2$ and $x \in \fun{U}_2\mo{M}$
  and $y \in \fun{U}_2\mo{M}'$.
\end{lem}

\noindent
  We omit the obvious proof. Observe that \eqref{eq:inst-le} simply says
  that two worlds are $\lan{L}_1$-logically equivalent if and only if
  they are $\lan{L}_2$-logically equivalent. Let us have a look at an example.

\begin{exm}\label{exm:inst-2}
  Let $\lan{L}_1 = \Lint$ and $\lan{L}_2 = \Lbiint$, both generated by the
  same set $\Prop$ of proposition letters.
  Since both can be interpreted in intuitionistic Kripke frames,
  there is an evident transformation $r : \mb{M}_2 \to \mb{M}_1$,
  namely the identity on the class of intuitionistic Kripke models.
  If we let $t : \Lint \to \Lbiint$ be the obvious translation,
  then clearly \eqref{eq:inst-th} is satisfied.
  However, the translation $t$ is not surjective.
  
  To overcome this, we can enrich $\Lint$ with an additional proposition letter
  $p_{\phi}$ for every formula $\phi \in \Lbiint$ that is not in $\Lint$.
  These can be interpreted by extending the valuation $V$ of an intuitionistic
  Kripke model $(X, \leq, V)$ via $V(p_{\phi}) = \llb \phi \rrb$,
  where the latter interpretation is given by the clauses in Definition
  \ref{def:interpr}. Denote this collection of additional proposition letters
  by $\Prop'$.
  Then clearly the translation $t : \Lint(\Prop) \to \Lbiint(\Prop)$
  extends to a surjective translation $t : \Lint(\Prop \cup \Prop') \to \Lbiint(\Prop)$.
  Moreover, we still have an obvious transformation of models and
  \eqref{eq:inst-th} is satisfied.
  It follows that the relation of $\Lbiint(\Prop)$-logical equivalence between
  two intuitionistic Kripke models coincides with
  $\Lint(\Prop \cup \Prop')$-logical equivalence.
\end{exm}

  More generally, if $\lan{L}_2$ freely extends $\lan{L}_1$ with one or more
  operators, then we can use this method to transfer properties of
  $\lan{L}_1$-logical equivalence to $\lan{L}_2$, achieving surjectivity
  by adding a proposition letter $p_{\phi}$ to $\lan{L}_1$ for
  every $\lan{L}_2$-formula that is not already in $\lan{L}_1$.
  
  We use this as follows: Suppose we know that logical equivalence between
  certain models for $\lan{L}_1$ is a bisimulation relation, and hence implies
  certain back-and-forth conditions. Then by the lemma the logical equivalence
  relation between models for $\lan{L}_2$ coincides with
  $\lan{L}_1$-logical equivalence, and therefore allows us to inherit the
  back-and-forth conditions.

\section{Bisimulations}\label{sec:non-modal}

  \noindent
  We begin this section by recalling the definition of bisimulation between Kripke models
  given in \cite{Pat97}, and prove a Hennessy-Milner result.
  We then dualise this to obtain a corresponding result for dual
  intuitionistic logic. Taken together, both results imply 
  the Hennessy-Milner property 
  for bi-intuitionistic logic.

\subsection{Bisimulations for Intuitionistic Logic}

\begin{defn}\label{def:i-bis}
  Let $\mo{M} = (X, \leq, V)$ and $\mo{M}' = (X', \leq', V')$ be two intuitionistic Kripke models.
  An \emph{intuitionistic bisimulation} or \emph{$\Lint$-bisimulation}
  between $\mo{M}$ and $\mo{M}'$ is a relation
  $B \subseteq X \times X'$ such that for all $(x, x') \in B$ we have:
  \begin{enumerate}[\qquad 1 \;]
    \renewcommand{\theenumi}{($B_{\arabic{enumi}}$)}
    \item \label{eq:bis1}
          For all $p \in \Prop$, $x \in V(p)$ iff $x' \in V'(p)$;
    \item \label{eq:bis2}
          If $x \leq y$ then there exists $y' \in X'$ such that $x' \leq' y'$ and $yBy'$;
    \item \label{eq:bis3}
          If $x' \leq' y'$ then there exists $y \in X$ such that $x \leq y$ and $yBy'$.
  \end{enumerate}
  Two states $x$ and $x'$ are called \emph{$\Lint$-bisimilar} if there is an $\Lint$-bisimulation
  linking them, notation: $x \bisim_{\Lint} x'$.
\end{defn}

  \noindent
  A straightforward inductive argument proves that bisimilar states
  satisfy the same formulae.

\begin{propn}\label{prop:int-bis-sound}
  If $x \bisim_{\Lint} x'$ then $x \logeq_{\Lint} x'$.
\end{propn}

  \noindent
  Furthermore, it is easy to see (but of no relevance for us in the
  sequel) that intuitionistic bisimulations are closed under composition,
  and the graph of a bounded morphism is an intuitionistic bisimulation.
  We prove a Hennessy-Milner property for image-compact models.

\begin{thm}\label{thm:hm-int}
  Let $x, x'$ be worlds in two image-compact models $\mo{M} = (X, \leq, V)$ and
  $\mo{M}' = (X', \leq', V')$. Then
  $$
    x \bisim_{\Lint} x' \iff x \logeq_{\Lint} x'.
  $$
\end{thm}
\begin{proof}
  Since we assume $\mo{M}$ and $\mo{M}'$ to be image-compact,
  they both carry a general model structure, i.e., we can
  find a collection $A$ of up-closed subsets of $(X, \leq)$ such that
  $(X, \leq, V, A)$ is a general model and ${\uparrow}_{\leq}x$ is
  compact in $\tau_A$ for all $x \in X$, and similarly for $\mo{M}'$.
  Suppose we have chosen such $A$ and $A'$.

  The direction from left to right is soundness of the notion of bisimulation
  (Proposition \ref{prop:int-bis-sound}).
  For the converse direction, we show that the relation of logical equivalence is
  a bisimulation between $\mo{M}$ and $\mo{M}'$.
  
  Clearly, if $x \logeq_{\Lint} x'$ we have $x \in V(p)$ iff $x' \in V'(p)$, so
  item \ref{eq:bis1} is satisfied.
  
  We now prove that \ref{eq:bis2} holds.
  Let $x \logeq_{\Lint} x'$ and $x \leq y$.
  Then we need to find $y' \in X'$ such that $x' \leq' y'$ and $y \logeq_{\Lint} y'$.
  Suppose towards a contradiction that such a $y'$ does not exist.
  Then for each $\leq'$-successor $z'$ of $x'$ we can either find a separating formula
  $\phi_{z'}$ such that $\mo{M}, y \Vdash \phi_{z'}$ and $\mo{M}', z' \not\Vdash \phi_{z'}$, or a
  separating formula $\psi_{z'}$ such that $\mo{M}, y \not\Vdash \psi_{z'}$ and
  $\mo{M}', z' \Vdash \psi_{z'}$. Pick such a separating formula for each $z'$. Let
  $\Phi$ be the collection of such formulae that are \emph{not} satisfied at $z'$,
  and $\Psi$ the collection of separating formulae that \emph{are} satisfied at $z'$.
  
  Since the interpretants of the formulae are clopen in the topology on $X'$ generated
  by $A' \cup -A'$, the collection
  $$
    \{ X \setminus \llb \phi \rrb^{\mo{M}'} \mid \phi \in \Phi \}
      \cup \{ \llb \psi \rrb^{\mo{M}'} \mid \psi \in \Psi \}
  $$
  is an open cover of ${\uparrow}_{\leq'}x'$. As the latter is assumed to be compact,
  we get finite subsets $\Phi' \subseteq \Phi$ and $\Psi' \subseteq \Psi$ such that
  $$
    \{ X \setminus \llb \phi \rrb^{\mo{M}'} \mid \phi \in \Phi' \}
      \cup \{ \llb \psi \rrb^{\mo{M}'} \mid \psi \in \Psi' \}
  $$
  covers ${\uparrow}_{\leq'}x'$. As a consequence, for every successor $z'$ of $x'$ there
  either exists a $\phi \in \Phi'$ such that $z' \not\Vdash \phi$, or a $\psi \in \Psi'$
  such that $z' \Vdash \psi$. Therefore,
  $$
    x' \Vdash \textstyle\bigwedge \Phi' \to \bigvee \Psi'.
  $$
  Since the disjunction and conjunction are taken over finite sets, this is a
  formula in $\Lint$.
  Furthermore, $y$ satisfies all $\phi \in \Phi'$ and none of the $\psi \in \Psi'$,
  and hence
  $$
    x \not\Vdash \textstyle\bigwedge \Phi' \to \bigvee \Psi'.
  $$
  This is a contradiction with the assumption that $x$ and $x'$ are logically equivalent.
  Therefore there must exist $y' \in X'$ which is logically equivalent to $y$ and
  satisfies $x' \leq' y'$.
  Item \ref{eq:bis3} is proven symmetrically.
\end{proof}

  \noindent
  Theorem \ref{thm:hm-int} does not give a strict characterisation of
  models where logical equivalence coincides with bisimilarity.
  This is witnessed by the following example, which gives a model that is \emph{not}
  image-compact while logical equivalence (between the model and itself) \emph{does}
  imply bisimilarity.
  
\begin{exm}
  Consider the intuitionistic Kripke frame consisting of the
  rational numbers ordered as usual. Let $\Prop = \{ p_q \mid q \in \mathbb{Q} \}$
  be a countable set of proposition letters and define a valuation
  $V : \Prop \to \fun{Up}(\mathbb{Q}, \leq)$ by
  $V(p_q) = \{ x \in \mathbb{Q} \mid q < x \}$.
  Then $\mo{Q} = (\mathbb{Q}, \leq, V)$ is an intuitionistic Kripke model.
  
  We claim that $\mo{Q}$ is not image-compact. To see this, let $A$ be any general frame structure
  such that $(\mathbb{Q}, \leq, A, V)$ is a general model.
  By definition $\llb p_q \rrb \in A \cup -A$ and $\mathbb{Q} \setminus \llb p_q \rrb \in A \cup -A$
  for all $p_q \in \Prop$. We note that
  ${\uparrow}_{\leq}0$ is covered by
  $$
    (\mathbb{Q} \setminus \llb p_0 \rrb)
      \cup \bigcup \{ \llb p_{\nicefrac{1}{n}} \rrb \mid n \in \mathbb{N} \}
  $$
  and clearly this cover does not have a finite subcover.
  However, the relation of logical equivalence between $\mo{Q}$ and itself
  is the identity, and hence is automatically an $\Lint$-bisimulation.
\end{exm}

  Also, it is not in general true that logical equivalence implies bisimilarity.
  In \cite[Proposition 27]{Pat97} the author gives an example of two
  intuitionistic Kripke models such that logical equivalence does not imply
  bisimilarity. (The notion of image-finiteness used in
  {\it loc.~\!cit.}~is not the usual one.)
  Alternatively, one can give a counterexample using ``porcupine models''
  similar to Example \ref{exm:biint1} below.

\subsection{Bisimulations for Dual- and Bi-Intuitionistic Logic}

\begin{defn}\label{def:d-bis}
  A \emph{dual-intuitionistic bisimulation} or \emph{$\Ldint$-bisimulation}
  between intuitionistic Kripke models $\mo{M} = (X, \leq, V)$ and $\mo{M}' = (X', \leq', V')$
  is a relation $B \subseteq X \times X'$ such that for all $(x, x') \in B$ we have:
  \begin{enumerate}[\qquad 1 \;]
    \renewcommand{\theenumi}{($B_{\arabic{enumi}}$)}
    \item For all $p \in \Prop$, $x \in V(p)$ iff $x' \in V'(p)$;
    \setcounter{enumi}{3}
    \item \label{eq:bis4}
          If $y \leq x$ then there exists $y' \in X'$ such that $y' \leq' x'$ and $yBy'$;
    \item \label{eq:bis5}
          If $y' \leq' x'$ then there exists $y \in X$ such that $y \leq x$ and $yBy'$.
  \end{enumerate}
  If moreover $B$ satisfies \ref{eq:bis2} and \ref{eq:bis3} (from Definition
  \ref{def:i-bis}) then we call $B$ a \emph{bi-intuitionistic bisimulation},
  or \emph{$\Lbiint$-bisimulation}.
  We define $\Ldint$-bisimilarity and $\Lbiint$-bisimilarity as usual, and write these as
  $x \bisim_{\Ldint} x'$ and $x \bisim_{\Lbiint} x'$.
\end{defn}

\begin{rem}\label{rem-bisim-Badia}
  \emph{Directed $\Lbiint$-bisimulations} \cite[Definition 4]{Bad16} between
  intuitionistic Kripke models are pairs $(Z_1, Z_2)$
  of simulations, i.e., pairs $(Z_1, Z_2)$ of two relations $Z_1
  \subseteq X \times X'$ and $Z_2 \subseteq X' \times X$ satisfying
  certain back-and-forth conditions. This is closely related to $\Lbiint$-bisimulation as just
  introduced: if $B$ is
  a $\Lbiint$-bisimulation then $(B, B^{-1})$ is a directed $\Lbiint$-bisimulation,
  and conversely if $(Z_1, Z_2)$ is a directed $\Lbiint$-bisimulation,
  then $Z_1 \cap Z_2^{-1}$ is a $\Lbiint$-bisimulation.

  Although not carried out in \emph{op.~\!cit.}, one could define
  $x$ and $x'$ to be \emph{directed $\Lbiint$-bisimilar}
  if there is a directed $\Lbiint$-bisimulation $(Z_1, Z_2)$ with 
  $(x, x') \in Z_1$ and $(x', x) \in Z_2$. Directed 
  $\Lbiint$-bisimilarity and $\Lbiint$-bisimilarity as defined in Definition
  \ref{def:d-bis} above are then easily seen to coincide.
\end{rem}

\begin{figure}[h!]
  \centering
    \begin{tikzcd}[arrows=-]
      y \arrow[r, dashed] & y' &
        y \arrow[r, dashed] & y' &
        x \arrow[r] & x' &
        x \arrow[r] & x' \\
      x \arrow[u, dashed, "\leq" left] \arrow[r, "B" below] & x' \arrow[u] &
        x \arrow[u, "\leq" left] \arrow[r, "B" below] & x' \arrow[u, dashed] &
        y \arrow[u, dashed, "\leq" left] \arrow[r, dashed, "B" below] & y' \arrow[u] &
        y \arrow[u, "\leq" left] \arrow[r, dashed, "B" below] & y' \arrow[u, dashed]
    \end{tikzcd}
  \caption{The zigs and zags of a $\Lbiint$-bisimulation.}
  \label{fig:zigzag1}
\end{figure}
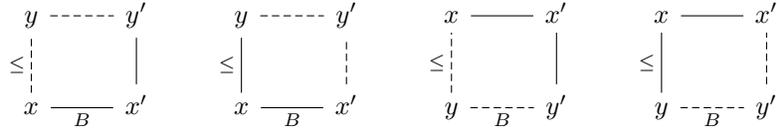

\begin{propn}\label{prop:b-d-bisim-sound}
  Let $(X, \leq, V)$ and $(X', \leq', V')$ be two intuitionistic Kripke models and
  $x \in X, x' \in X'$.
  Then $x \bisim_{\Ldint} x'$ implies $x \logeq_{\Ldint}$ and
  $x \bisim_{\Lbiint} x'$ implies $x \logeq_{\Lbiint} x'$.
\end{propn}

  \noindent
  The following lemma allows us to view an $\Ldint$-bisimulation between two
  models $\mo{M}$ and $\mo{M}'$ as an $\Lint$-bisimulation between the corresponding
  dual models.

\begin{lem}\label{lem:d-i-bisim}
  Let $\mo{M} = (X, \leq, V)$ and $\mo{M}' = (X', \leq', V')$
  be two intuitionistic Kripke models. Then
  $B \subseteq X \times X'$ is a $\Ldint$-bisimulation between $\mo{M}$ and $\mo{M}'$
  if and only if it is an $\Lint$-bisimulation between $\mo{M}^{\partial}$
  and $(\mo{M}')^{\partial}$.
\end{lem}

  \noindent
  Using this lemma we can convert the result from Theorem \ref{thm:hm-int}
  to a Hennessy-Milner theorem for dual-intuitionistic logic.

\begin{thm}\label{thm:hm-d-int}
  Let $x, x'$ be worlds in two pre-image-compact intuitionistic Kripke models
  $\mo{M} = (X, \leq, V)$ and $\mo{M}' = (X', \leq', V')$. Then
  $$
    x \bisim_{\Ldint} x' \iff x \logeq_{\Ldint} x'.
  $$
\end{thm}
\begin{proof}
  Let $B$ be the relation of logical equivalence between $\mo{M}$ and $\mo{M}$.
  We show that it is an $\Ldint$-bisimulation.
  Alternatively, it suffices to show that it is an $\Lint$-bisimulation
  between $\mo{M}^{\partial}$ and $(\mo{M}')^{\partial}$.

  By Lemma \ref{lem:sem-trans}, two states $x, x'$ in $\mo{M}^{\partial}$ 
  and $(\mo{M}')^{\partial}$ satisfy the same $\Lint$-formulae
  if and only if they satisfy the same $\Ldint$-formulae
  in $\mo{M}$ and $\mo{M}'$. Therefore the relation $B$ coincides with
  logical equivalence between $\mo{M}^{\partial}$ and $(\mo{M}')^{\partial}$.
  Furthermore, $\mo{M}^{\partial}$ and $(\mo{M}')^{\partial}$
  are image-compact because $\mo{M}$ and $\mo{M}'$ are pre-image-compact.
  So it follows from Theorem
  \ref{thm:hm-int} that $B$ is an $\Lint$-bisimulation between
  $\mo{M}^{\partial}$ and $(\mo{M}')^{\partial}$,
  hence an $\Ldint$-bisimulation between $\mo{M}$ and $\mo{M}'$.
\end{proof}

  \noindent
  Combining Lemma \ref{lem:inst} and Theorems \ref{thm:hm-int} and \ref{thm:hm-d-int}
  yields:
  
\begin{thm}\label{thm:hm-bi-int}
  Let $x, x'$ be worlds in two intuitionistic Kripke
  models $\mo{M} = (X, \leq, V)$ and
  $\mo{M}' = (X', \leq', V')$ that are both image-compact and pre-image-compact. Then
  $$
    x \bisim_{\Lbiint} x' \iff x \logeq_{\Lbiint} x'.
  $$
\end{thm}
\begin{proof}
  The direction from left to right follows is Proposition \ref{prop:b-d-bisim-sound}.
  For the converse, we will show that the relation $B$ of logical equivalence between them is
  a bisimulation.
  \ref{eq:bis1} follows immediately from the fact
  that $B$ is logical equivalence.
  
  To show that \ref{eq:bis2} and \ref{eq:bis3} hold, we use Lemma \ref{lem:inst}.
  Let $\Prop'$ be defined as in Example \ref{exm:inst-2} and extend the valuations
  $V$ and $V'$ of $\mf{M}$ and $\mf{M}'$ to $\hat{V}$ and $\hat{V}'$
  by setting $\hat{V}(p_{\phi}) = \llb \phi \rrb^{\mf{M}}$, and similar for $\hat{V}'$.
  Then as a consequence of Lemma \ref{lem:inst}, $B$ coincides with
  $\Lint$-logical equivalence between $(X, \leq, \hat{V})$ and $(X', \leq', \hat{V}')$.
  Furthermore, these new models are image-compact, and therefore
  properties \ref{eq:bis2} and \ref{eq:bis3} follow from Theorem \ref{thm:hm-int}.
  
  One can similarly obtain \ref{eq:bis4} and \ref{eq:bis5} from Theorem \ref{thm:hm-d-int}.
\end{proof}

  \noindent
  We complete this section with a detailed example
  showing that logical equivalence for bi-intuitionistic formulae does not
  in general imply $\Lbiint$-bisimilarity.

\begin{figure}
\centering
\begin{tikzpicture}[shorten <=4pt, shorten >=4pt]
  \draw (-.7,.7) node{$\mf{W}$};
  \draw[fill=black] (0,0) circle(1.6pt)
                        node[anchor=south east]{\footnotesize{$q$}};
  \draw[fill=black] (1,0) circle(1.6pt)
                        node[anchor=south west]{\footnotesize{$p_0$}};
  \draw[fill=black] (-30:1) circle(1.6pt)
                        node[anchor=south west]{\footnotesize{$p_1$}}
                    (-30:2) circle(1.6pt)
                        node[anchor=south west]{\footnotesize{$p_0$}}
                        node[anchor=north, rotate=5]{\scriptsize{$(2,1)$}};
  \draw[fill=black] (-60:1) circle(1.6pt)
                        node[anchor=west]{\footnotesize{$p_2$}}
                    (-60:2) circle(1.6pt)
                        node[anchor=west]{\footnotesize{$p_1$}}
                        node[anchor=east, rotate=15]{\scriptsize{$(3, 1)$}}
                    (-60:3) circle(1.6pt)
                        node[anchor=west]{\footnotesize{$p_0$}}
                        node[anchor=east, rotate=15]{\scriptsize{$(3, 0)$}};
  \draw[fill=black] (-100:1.5) circle(1.6pt)
                        node[anchor=east]{\footnotesize{$p_3$}}
                    (-100:2.5) circle(1.6pt)
                        node[anchor=east]{\footnotesize{$p_2$}}
                        node[anchor=west, rotate=-3]{\scriptsize{$(\infty, 2)$}}
                    (-100:3.5) circle(1.6pt)
                        node[anchor=east]{\footnotesize{$p_1$}}
                        node[anchor=west, rotate=-3]{\scriptsize{$(\infty, 1)$}}
                    (-100:4.5) circle(1.6pt)
                        node[anchor=east]{\footnotesize{$p_0$}}
                        node[anchor=west, rotate=-3]{\scriptsize{$(\infty, 0)$}};
  \draw[latex-] (0,0) -- (1,0);
  \draw[latex-] (0,0) -- (-30:1);
  \draw[latex-] (-30:1) -- (-30:2);
  \draw[latex-] (0,0) -- (-60:1);
  \draw[latex-] (-60:1) -- (-60:2);
  \draw[latex-] (-60:2) -- (-60:3);
  \draw[dashed] (0,0) -- (-100:1.5);
  \draw[latex-] (-100:1.5) -- (-100:2.5);
  \draw[latex-] (-100:2.5) -- (-100:3.5);
  \draw[latex-] (-100:3.5) -- (-100:4.5);
  \draw[thick, dotted] (-67:1.7) arc(-67:-94:1.7);
\end{tikzpicture}
\qquad\qquad
\begin{tikzpicture}[shorten <=4pt, shorten >=4pt]
  \draw (-.7,.7) node{$\mf{W}'$};
  \draw[fill=black] (0,0) circle(2pt)
                        node[anchor=south east]{\footnotesize{$q$}};
  \draw[fill=black] (1,0) circle(2pt)
                        node[anchor=south west]{\footnotesize{$p_0$}};
  \draw[fill=black] (-30:1) circle(2pt)
                        node[anchor=south west]{\footnotesize{$p_1$}}
                    (-30:2) circle(2pt)
                        node[anchor=south west]{\footnotesize{$p_0$}}
                        node[anchor=north, rotate=5]{\scriptsize{$(2,1)$}};
  \draw[fill=black] (-60:1) circle(2pt)
                        node[anchor=west]{\footnotesize{$p_2$}}
                    (-60:2) circle(2pt)
                        node[anchor=west]{\footnotesize{$p_1$}}
                    (-60:3) circle(2pt)
                        node[anchor=south west]{\footnotesize{$p_0$}}
                        node[anchor=north, rotate=5]{\scriptsize{$(3,1)$}};
  \draw[fill=black] (-80:1) circle(2pt)
                        node[anchor=east]{\footnotesize{$p_3$}}
                    (-80:2) circle(2pt)
                        node[anchor=west]{\footnotesize{$p_2$}}
                        node[anchor=east, rotate=5]{\scriptsize{$(4,2)$}}
                    (-80:3) circle(2pt)
                        node[anchor=west]{\footnotesize{$p_1$}}
                        node[anchor=east, rotate=5]{\scriptsize{$(4,1)$}}
                    (-80:4) circle(2pt)
                        node[anchor=west]{\footnotesize{$p_0$}}
                        node[anchor=east, rotate=5]{\scriptsize{$(4,0)$}};
  \draw[latex-] (0,0) -- (1,0);
  \draw[latex-] (0,0) -- (-30:1);
  \draw[latex-] (-30:1) -- (-30:2);
  \draw[latex-] (0,0) -- (-60:1);
  \draw[latex-] (-60:1) -- (-60:2);
  \draw[latex-] (-60:2) -- (-60:3);
  \draw[latex-] (0,0) -- (-80:1);
  \draw[latex-] (-80:1) -- (-80:2);
  \draw[latex-] (-80:2) -- (-80:3);
  \draw[latex-] (-80:3) -- (-80:4);
  \draw[thick, dotted] (-100:1.7) arc(-100:-119:2);
  \draw[opacity=0] (-100:4.5) circle(2pt) node[anchor=east]{\footnotesize{$p_0$}}
                        node[anchor=west, rotate=-3]{\scriptsize{$(\infty, 0)$}};
\end{tikzpicture}
\caption{The figures depicts the models $\mf{W}$ and $\mf{W}'$ from Example \ref{exm:biint1}.
The coordinates indicate the names of some of the states. The $p_i$ denote the lowest occurrence of a proposition letter in each branch of the models. That is, if $p_i$ is true in some state, then it is also true in all states above.}
\label{fig-exm-fin-comp}
\end{figure}
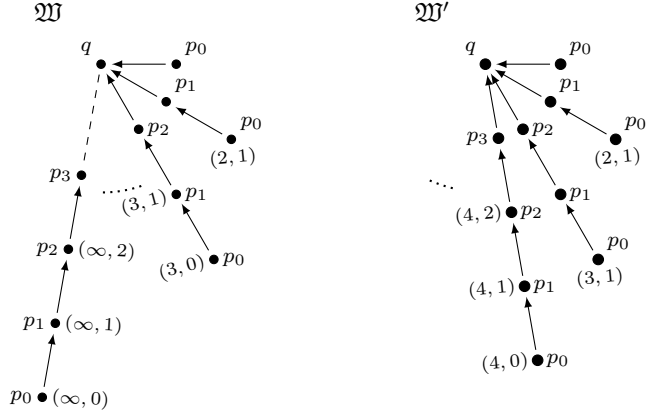

\begin{exm}\label{exm:biint1}
  Let $W = \{ (n, k) \in (\mathbb{N}
           \cup \{ \infty \}) \times \mathbb{N} \mid k < n \}
           \cup \{ x \}$
  and define an order $\preq$ by: $(n, k) \preq x$ for all $(n, k) \in W$ and 
  $(n, k) \preq (m, \ell)$ iff $n = m$ and $k \leq \ell$.
  For $\Prop = \{ p_i \mid i \in \mb{N} \} \cup \{ q \}$ define the valuation 
  $V$ by $V(q) = \{ x \}$ and $V(p_i) = \{ (n, k) \in W \mid i \leq k \} \cup \{ x \}$. 
  Then the triple $(W, \preq, V)$ is a Kripke model.

  Let $\mf{W}' = (W', \preq', V')$ be the submodel of $\mf{W}$ with underlying 
  set $W' = \{ (n', k') \in \mathbb{N} \times \mathbb{N} \mid k' < n' \} \cup \{ x' \}$. 
  Note that $\mf{W}'$ does not have an infinite branch.
  (We use primes to distinguish the two models.)
  See Figure \ref{fig-exm-fin-comp} for pictorial presentations of the two models.

  We claim that $x$ and $x'$ are logically equivalent but not bisimilar.
  Suppose towards a contradiction that there exists a bisimulation $B$ linking 
  $x$ and $x'$. Since $(\infty, 0) \preq x$ in $W$ there must be some $y' \in W'$ 
  such that $(\infty, 0)By'$ and $y' \preq' x'$. 
  Then $y'$ cannot be $x'$, because $W, (\infty, 0) \not\Vdash q$, 
  hence $W', y' \not\Vdash q$, whereas $W', x' \Vdash q$. 
  So $y'$ is of the form $(n', k')$ for some $n', k' \in \mb{N}$ with $k' < n'$. 
  But then $W', (n', k') \Vdash p_{n'+1} \to q$, 
  while $W, (\infty, 0) \not\Vdash p_{n'+1} \to q$. Therefore $(\infty, 0)$ and 
  $(n', k')$ are not logically equivalent, hence by
  Proposition \ref{prop:b-d-bisim-sound} they cannot be bisimilar.
  This contradicts the assumption that there exists a bisimulation
  $B$ linking $x$ and $x'$, thus $x$ and $x'$ are not bisimilar.

  Next we show that $x \in W$ and $x' \in W'$ are logically equivalent.
  For $m \in \mb{N}$, let $\Prop_m = \{ p_i \mid i \in \mb{N}, i \leq m \} \cup \{ q \}$. 
  Then $\Lint(\Prop) = \bigcup_{m \in \mb{N}} \Lint(\Prop_m)$.
  Define $B_m \subseteq W \times W'$ by
  \begin{align*}
    B_m = \{ (x, x') \} \cup \big\{ \big((n,k), (n',k')\big) \mid
      \text{either } &[ n = n' \text{ and } k = k' ] \\
      \text{or } &[k, k' \geq m ] \\
      \text{or } &[n, n' > m \text{ and } k = k' < m ] 
      \big\}.
  \end{align*}
  It can be shown by induction that whenever $(z, z') \in B_m$, 
  we have $W, z \Vdash \phi$ iff $W', z' \Vdash \phi$ for all $\phi \in \mb{L}(\Prop_m)$.
  It follows that $x$ and $x'$ are logically equivalent because $(x, x') \in B_m$ 
  for all $m \in \mb{N}$. As we have already established that $x$ and
  $x'$ are not bisimilar, we conclude that logical equivalence
  cannot imply bisimilarity in general.
\end{exm}

\section{Modal Bi-/Dual-/Intuitionistic Logics}\label{sec:modal}

  \noindent
  In this section we enrich the logics from Section \ref{sec:non-modal} with
  (several copies of) the unary modal operators $\Box$ and $\Diamond$.
  Following \cite{BozDos84}, we shall treat $\Box$ and $\Diamond$ as two different
  modalities that a priori are not related via axioms.
  Semantically, $\Box$ and $\Diamond$ are interpreted via distinct
  relations, so that 
  boxes and diamonds do not necessarily come in pairs.
  For $\lan{L} \in \{ \Lint, \Ldint, \Lbiint \}$,
  we write $\lan{L}_{n,m}$ for the languages that arises from adjoining
  $\lan{L}$ with boxes $\Box_1, \ldots, \Box_n$ diamonds $\Diamond_1, \ldots, \Diamond_m$.
  In the special case where $n = 1$ and $m = 0$ we write
  $\lan{L}_{\Box} := \lan{L}_{1,0}$, and similarly we sometimes use
  $\lan{L}_{\Diamond} := \lan{L}_{0,1}$
  and $\lan{L}_{\Box\Diamond} := \lan{L}_{1,1}$.

  Since we do not assume any axioms relating boxes and diamonds, each
  modality is interpreted via its own relation in the same way as in classical
  modal logic.
  As such, a model for $\lan{L}_{n,m}$ is an intuitionistic Kripke model
  with an additional relation $R_i$ for each box and $S_j$ for each diamond,
  satisfying certain coherence conditions with respect to the order $\leq$
  to ensure that the interpretation of every formula is an upset.
  This approach resembles that of $H\Box$- and $H\Diamond$-models introduced
  in \cite{BozDos84}.

  The main objective of this section is to obtain a Hennessy-Milner type theorem for the
  modal bi-intuitionistic logic $\lan{L}_{n,m}$ interpreted in the models sketched above.
  We shall prove intermediate results for $\Lint_{\Box} = \Lint_{1,0}$ and 
  $\Ldint_{\Diamond} = \Ldint_{0,1}$,
  which we then combine for the desired result using Lemma \ref{lem:inst}.

\subsection{Semantics for Modal Bi-/Dual-/Intuitionistic Logics}
  
  \noindent
  If $Z$ and $Z'$ are two relations on a set $X$, then we denote by
  $Z \circ Z'$ the relation
  $\{ (x, y) \in X \times X \mid \exists u \in X \text{ s.t. } xZu \text{ and } uZ'y \}$.
  
\begin{defn}\label{def:modal-model}
  A \emph{(modal) $\lan{L}_{n,m}$-frame} is a tuple
  $(X, \leq, R_1, \ldots, R_n, S_1, \ldots, S_m)$ that consists of an intuitionistic
  Kripke frame $(X, \leq)$ and relations $R_i, S_j \subseteq X \times X$ satisfying
  $$
    ({\leq} \circ R_i) \subseteq (R_i \circ {\leq}),
    \qquad({\geq} \circ S_j) \subseteq (S_j \circ {\geq}).
  $$
  It is called \emph{strictly condensed} if
  $
    ({\leq} \circ R_i \circ {\leq}) \subseteq R_i
  $
  and
  $
    ({\geq} \circ S_j \circ {\geq}) \subseteq S_j
  $
  for all $i \in \{ 1, \ldots, n \}$ and $j \in \{ 1, \ldots, m \}$.
  The corresponding notion of an \emph{$\lan{L}_{n,m}$-model} arises from adding a valuation.
\end{defn}

  Note that, since $\leq$ is reflexive, an $\lan{L}_{n,m}$-frame is strictly condensed
  if and only if $({\leq} \circ R_i \circ {\leq}) = R_i$
  and $({\geq} \circ S_j \circ {\geq}) = S_j$ for all $i$ and $j$.

  These models can be used to interpret modal extensions of
  $\Lint, \Ldint$ and $\Lbiint$ with $n$ boxes and $m$ diamonds.
  The logical connectives from $\lan{L}$ are interpreted in the underlying
  intuitionistic Kripke model $(X, \leq, V)$ as usual and,
  as stated, the interpretations of $\Box_i$ and $\Diamond_j$ are defined as in
  classical modal logic, via the relations $R_i$ and $S_j$.
  That is,
  \begin{align*}
    \mo{M}, x \Vdash \Box_i\phi &\iff \text{for all } y \in X, \; xR_iy \text{ implies } \mo{M}, y \Vdash \phi \\
    \mo{M}, x \Vdash \Diamond_j\phi &\iff \mo{M}, y \Vdash \phi
    \text{ for some $y$ with } xS_jy.
  \end{align*}
  We write $x \logeq_{\lan{L}_{n,m}} x'$ if two states
  satisfy precisely the same $\lan{L}_{n,m}$-formulae.

  We shall sometimes write
  $(X, \leq, (R_i), (S_j), V)$ for a modal $\lan{L}_{n,m}$-model.
  Besides, we remark that (strictly condensed) $\lan{L}_{\Box}$-models are precisely
  (strictly condensed) $H\Box$-models from \cite{BozDos84},
  and (strictly condensed) $\lan{L}_{\Diamond}$-models can be found in {\it op.~\!cit.}~under
  the name of (strictly condensed) $H\Diamond$-models.
  We have the following notion of bisimulation for these models:
  
\begin{defn}\label{def:modal-bis}
  Let $\mo{M} = (X, \leq, (R_i), (S_j), V)$ and $\mo{M}' = (X', \leq', (R_i'), (S_j'), V')$
  be two modal $\lan{L}_{n,m}$-models and $B \subseteq X \times X'$ a relation.
  We call $B$ a $\Box_i$-zigzag if for all $(x, x') \in B$ the following conditions hold:
  \begin{description}
    \item[($\Box_i$-zig)] If $xR_iy$ then there exists $y' \in X'$ such that
          $x'R_i'y'$ and $yBy'$;
    \item[($\Box_i$-zag)] If $x'R'_iy'$ then there exists $y \in X$ such that
          $xR_iy$ and $yBy'$;
  \end{description}
  We call $B$ a $\Diamond_j$-zigzag if the same conditions hold for $S_j$
  instead of $R_i$.

  An \emph{$\lan{L}_{n,m}$-bisimulation} between $\mo{M}$ and $\mo{M}'$ is a
  relation $B \subseteq X \times X'$ which is an $\lan{L}$-bisimulation
  between the underlying intuitionistic Kripke models and which is a
  $\Box_i$-zigzag and $\Diamond_j$-zigzag for all $i \in \{ 1, \ldots, n \}$
  and $j \in \{ 1, \ldots, m \}$.
\end{defn}

  We remark that one can quotient with bisimilarity:

\begin{rem}
  Let $\mf{M} = (X, \leq, (R_i), (S_j), V)$ be an $\lan{L}_{n,m}$-model.
  It is easy to see that the collection of $\lan{L}_{n,m}$-bisimulation on a model
  is closed under all unions.
  Therefore, the relation $B$ of bisimilarity on $\mf{M}$ is again a bisimulation.
  Moreover, $B$ is an equivalence relation: it is reflexive because the identity
  on $X$ is a bisimulation, symmetric because the inverse of a bisimulation on $X$
  is again a bisimulation, and transitive because bisimulations are closed under
  composition.
  
  Let $X_B$ denote the quotient of $X$ with the equivalence relation $X$
  and write $\bar{x} \in X_B$ for the equivalence class of $x \in X$.
  For each of the relations $Z$ on $X$, define a relation $Z_B$ on $X/B$ via
  $\bar{x}Z_B\bar{y}$ if there are $x' \in \bar{x}$ and $y' \in \bar{y}$
  such that $x'Zy'$. Finally, for $p \in \Prop$ let $V_B(p) = \{ \bar{x} \mid x \in V(p) \}$.
  Then it follows from a straightforward verification that the tuple
  $$
    \mf{M}_B = (X_B, \leq_B, ((R_B)_i), ((S_B)_i), V)
  $$
  is an $\lan{L}_{n,m}$-model
  and the graph of the quotient map $q : X \to X_B$ is a bisimulation
  between $\mf{M}$ and $\mf{M}_B$.
  Consequently, if $\mf{M}$ is in a Hennesy-Milner class, then we can that the
  quotient with respect to logical equivalence.
\end{rem}

\begin{rem}
  When equipped with a suitable notion of (bounded) morphism, the collection of 
  $\lan{L}_{n,m}$-frames forms a category.
  This category is isomorphic to a category of \emph{dialgebras} \cite{GroPat20},
  and the language $\lan{L}_{n,m}$ arises as a \emph{dialgebraic logic}.
  Interestingly, on the level of frames,
  the bisimulations defined in Definition \ref{def:modal-bis}
  correspond precisely to \emph{dialgebra bisimulations} (or \emph{cospans}) in the 
  category of $\lan{L}_{n,m}$-frames.
\end{rem}

  \noindent
  A straightforward inductive proof yields:

\begin{propn}\label{prop:adeq-modal}
  Let $x$ and $x'$ be two states in $\lan{L}_{n,m}$-models $\mo{M}$ and $\mo{M}'$.
  Then $x \bisim_{\lan{L}_{n,m}} x'$ implies $x \logeq_{\lan{L}_{n,m}} x'$.
\end{propn}

  \noindent
  In order to get a suitable notion of (pre-)image-compactness for the relations
  $R_i, S_j$ we extend the notion of a general frame to this modal setting.

\begin{defn}\label{def:modal-general}
  A \emph{general $\lan{L}_{n,m}$-frame} consists of a modal $\lan{L}_{n,m}$-frame
  $(X, \leq, R_1, \ldots, R_n, S_1, \ldots, S_m)$ and
  a collection $A \subseteq \fun{Up}(X, \leq)$ such that $(X, \leq, A)$ is a
  general $\lan{L}$-frame and $A$ is closed under:
  \begin{align*}
    \cbox_i &: \fun{Up}(X, \leq) \to \fun{Up}(X, \leq)
             : a \mapsto \{ x \in X \mid R_i[x] \subseteq a \} \\
    \cdiamond_j &: \fun{Up}(X, \leq) \to \fun{Up}(X, \leq)
                 : a \mapsto \{ x \in X \mid xS_jy \text{ for some } y \in a \}
  \end{align*}
  for all $i \in \{ 1, \ldots, n \}$ and $j \in \{ 1, \ldots, m \}$.
  The corresponding notion of a \emph{general $\lan{L}_{n,m}$-model}
  arises from adjoining such a frame
  with an \emph{admissible} valuation, i.e., a map $V : \Prop \to A$.

  A relation $R_i$ in an $\lan{L}_{n,m}$-model $(X, \leq, (R_i), (S_j), V)$
  is called \emph{(pre-)image-compact}
  if there exists $A \subseteq \fun{Up}(X, \leq)$
  such that $(X, \leq, (R_i), (S_j), A, V)$ is a general $\lan{L}_{n,m}$-model
  and $R_i[x] = \{ y \in X \mid xR_iy \}$ (resp. $R^{-1}_i[x] = \{ y \in X \mid yR_ix \}$)
  is compact in $\tau_A$ for every $x \in X$. 
  We similarly define (pre-)image-compactness for $S_j$.
\end{defn}

\begin{rem}
  The definition of (pre-)image-compactness crucially depends on
  the underlying base logic. In particular, we never speak about an
  image compact relation in an intuitionistic Kripke frame: we speak
  about an image compact relation in a $\Lint$-, $\Ldint$- or
  $\Lbiint$-frame. For a relation to qualify as image compact, we need to
  exhibit a system $A$ of admissible subsets that is \emph{closed
  under the operations of the base logic}. That is, a choice of
  admissibles may exhibit a relation as image compact in an
  $\Lint$-frame, but there may be no choice of admissibles $A'$ that
  exhibits the same relation as image-compact in an $\Lbiint$-frame,
  for example, if $A$ is not closed under $\undbito$.
  This subtlety is caused by the fact that we treat three base logics
  simultaneously.
\end{rem}

  For our Hennessy-Milner type results, we need to restrict to the
  strictly condensed models. Although this may seem like a harsh restriction,
  in fact every $\lan{L}_{n,m}$-model can be turned into a strictly condensed one without changing
  the interpretation of formulae, by merely readjusting the relations $R_i$ and $S_j$.
  We explicitly give this construction for $\lan{L}_{\Box}$-models and
  leave the general case to the reader.
  
\begin{propn}\label{prop:strictify}
  Let $\mo{M} = (X, \leq, R, V)$ be an $\lan{L}_{\Box}$-model and set
  $R^+ := ({R} \circ {\leq})$.
  Then $\mo{M}^+ = (X, \leq, R^+, V)$ is strictly condensed, and for all
  $x \in X$ and $\phi \in \lan{L}_{\Box}$ we have $\mo{M}, x \Vdash \phi$
  iff $\mo{M}^+, x \Vdash \phi$.
\end{propn}
\begin{proof}
  To see that $\mo{M}^+$ is strictly condensed, observe that reflexivity
  and transitivity of $\leq$ imply
  $({\leq} \circ R^+ \circ {\leq})
    = ({\leq} \circ R \circ {\leq})
    \subseteq 
    (R \circ {\leq})
    = R^+$.
  The preservation of truth can be proved by induction on the structure of the formula $\phi$.
  All cases are trivial except the modal case. For this, we have
  \begin{align*}
    \mo{M}, x \Vdash \Box\phi
      &\iff \text{for all $y \in X$, } xRy \text{ implies } \mo{M}, y \Vdash \phi \\
      &\iff \text{for all $y \in X$, } x(R \circ {\leq})y \text{ implies } \mo{M}, y \Vdash \phi \\
      &\iff \text{for all $y \in X$, } xR^+y \text{ implies } \mo{M}^+, y \Vdash \phi \\
      &\iff \mo{M}^+, x \Vdash \Box\phi.
  \end{align*}
  The second ``iff'' holds by the fact that truth-sets of formulae are up-closed in
  $(X, \leq)$, the third one by induction.
\end{proof}

  \noindent
  An example of this procedure is depicted in Figure \ref{fig:exm-M-plus} below.
  It is not in general true that either the identity or the 
  relation of logical equivalence between $\mo{M}$ and $\mo{M}^+$ is a
  $\lan{L}_{\Box}$-bisimulation, as is witnessed by the following example.

\begin{exm}\label{exm:M-plus}
  Let $X = \{ x, y, z \}$ be ordered by the pre-order generated by $y \leq z$
  and let $R = \{ (x, y) \} \subseteq X \times X$. Then $(X, \leq, R)$ is an
  $\lan{L}_{\Box}$-frame. Equip this with the valuation $V : \{ p, q \} \to \fun{Up}(X, \leq)$
  given by $V(p) = \{ y, z \}$ and $V(q) = \{ z \}$. Then
  $\mo{M} = (X, \leq, R, V)$ is the $\lan{L}_{\Box}$-model depicted in Figure \ref{fig:exm-M-plus}.
  The strictly condensed $\lan{L}_{\Box}$-model $\mo{M}^+$ is obtained by changing $R$
  to $R^+ = (R \circ {\leq}) = \{ (x, y), (x, z) \}$.
  
  The the relation of logical equivalence between $\mo{M}$ and $\mo{M}^+$
  is simply the identity relation on $X$.
  It is easy to see that this is \emph{not} an $\lan{L}_{\Box}$-bisimulation: in $\mo{M}^+$
  there is an $R$-transition from $x$ to $z$. The only state in $\mo{M}$ that
  is logically equivalent to $z$ is $z$. But there is no $R_{\Box}$-transition
  from $x$ to $z$ in $\mo{M}$.
  So there can be no $\lan{L}_{\Box}$-bisimulation linking $x$ and $x'$.
\end{exm}

\begin{figure}[h!]
  \centering
    \begin{tikzcd}[column sep=.5em, row sep=1.4em, arrows=-latex]
      \mo{M}
        &
        & z
        & [4em] \mo{M}^+
        & 
        & z \\
        & y \arrow[ru, "\leq"]
        &
        &
        & y \arrow[ru, "\leq"]
        & \\
        & x \arrow[u, "R"]
        &
        &
        & x \arrow[u, "R^+"]
            \arrow[ruu, bend right=15, "R^+" swap]
        &
    \end{tikzcd}
  \caption{An $\lan{L}_{\Box}$-model and its condensed version.}
  \label{fig:exm-M-plus}
\end{figure}
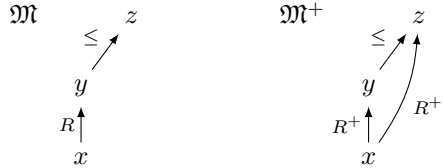

\subsection{Hennessy-Milner Property for Some Modal Intuitionistic Logics}
  
  \noindent
  We now restrict our attention to $\lan{Int_{\Box}}$ and extend the Hennessy-Milner result
  from Theorem \ref{thm:hm-int} to the setting of $\Lint_{\Box}$ interpreted in
  strictly condensed $\Lint_{\Box}$-models.

\begin{thm}\label{thm:hm-int-box}
  Let $\mo{M} = (X, \leq, R, V)$ and $\mo{M}' = (X', \leq', R', V')$
  be two strictly condensed $\Lint_{\Box}$-models such that $\leq, \leq', R$ and
  $R'$ are image-compact.
  Then for all $x \in X$ and $x' \in X'$ we have
  $$
    x \bisim_{\Lint_{\Box}} x' \iff x \logeq_{\Lint_{\Box}} x'.
  $$
\end{thm}
\begin{proof}
  The direction from left to right follows from Proposition \ref{prop:adeq-modal}.
  For the converse, we let $B$ be logical equivalence and we show that it 
  is a $\Lint_{\Box}$-bisimulation.
  It follows from Lemma \ref{lem:inst} and Theorem \ref{thm:hm-int} that $B$ 
  is an $\Lint$-bisimulation,
  so it remains to show that ($\Box$-zig) and ($\Box$-zag) hold.
  
  Let $xBx'$ and $xRy$ and suppose towards a contradiction that
  there is no $R'$-successor
  $y'$ of $x'$ which is logically equivalent to $y$. Then for each such $y'$ we can
  find a separating formula. As in Theorem \ref{thm:hm-int}, using compactness,
  we get two finite sets
  $\Phi'$ and $\Psi'$ such that $y$ satisfies every formula in $\Phi'$ and none in $\Psi'$,
  and such that for every $y'$ with $x'R'y'$ there either exists $\phi \in \Phi'$ such that
  $\mo{M}', y' \not\Vdash \phi$, or $\psi \in \Psi'$ such that $\mo{M}', y' \Vdash \psi$.
  
  Let $y'$ be an $R'$-successor of $x'$, then $y' \leq' z'$ implies
  $x'R'z'$, because $\mf{M}'$ is assumed to be strictly condensed.
  As a consequence $\mo{M}', y' \Vdash \bigwedge \Phi' \to \bigvee \Psi'$.
  Since this holds for any $y'$ with $x'R'y'$, we have
  $$
    \mo{M}', x' \Vdash \textstyle\Box(\bigwedge \Phi' \to \bigvee \Psi').
  $$
  Furthermore, by construction $\mo{M}, y \not\Vdash \bigwedge \Phi' \to \bigvee \Psi'$, so
  $$
    \mo{M}, x \not\Vdash \textstyle\Box(\bigwedge \Phi' \to \bigvee \Psi').
  $$
  This contradicts the assumption that $x$ and $x'$ be logically equivalent.
  Therefore we conclude that there must exist a $y' \in X'$ which is logically equivalent to
  $y$ and satisfies $x'R'y'$. Thus ($\Box$-zig) is satisfied.
  A symmetric argument shows that ($\Box$-zag) is satisfied as well.
\end{proof}

  \noindent
  The next example shows that
  a simple adaptation of ``porcupine models'' exhibits that logical equivalence
  does not in general imply $\Lint_{\Box}$-bisimilarity. Note also that in this example
  both $\leq$ and $\leq'$ are image-finite and pre-image-finite.

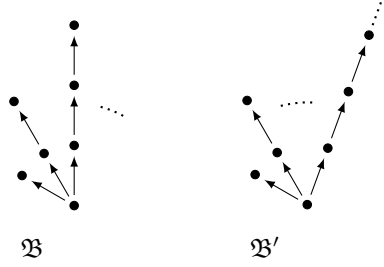
\begin{figure}[h!]
\centering
  \begin{tikzpicture}[scale=.8,shorten <=4pt, shorten >=4pt]
    \draw[fill=black] (0,0) circle(2pt);
    \draw[fill=black] (150:1) circle(2pt);
    \draw[fill=black] (120:1) circle(2pt)
                      (120:2) circle(2pt);
    \draw[fill=black] (90:1) circle(2pt)
                      (90:2) circle(2pt)
                      (90:3) circle(2pt);
    \draw[-latex] (0,0) -- (150:1);
    \draw[-latex] (0,0) -- (120:1);
    \draw[-latex] (120:1) -- (120:2);
    \draw[-latex] (0,0) -- (90:1);
    \draw[-latex] (90:1) -- (90:2);
    \draw[-latex] (90:2) -- (90:3);
    \draw[dotted, thick, shorten <=0pt, shorten >=0pt] (75:1.7) arc(75:60:1.7);
    \draw (-.7,-.7) node{$\mo{B}$};
  \end{tikzpicture}
  \qquad\qquad
  \begin{tikzpicture}[scale=.8,shorten <=4pt, shorten >=4pt]
    \draw[fill=black] (0,0) circle(2pt);
    \draw[fill=black] (150:1) circle(2pt);
    \draw[fill=black] (120:1) circle(2pt)
                      (120:2) circle(2pt);
    \draw[fill=black] (70:1) circle(2pt)
                      (70:2) circle(2pt)
                      (70:3) circle(2pt);
    \draw[-latex] (0,0) -- (150:1);
    \draw[-latex] (0,0) -- (120:1);
    \draw[-latex] (120:1) -- (120:2);
    \draw[-latex] (0,0) -- (70:1);
    \draw[-latex] (70:1) -- (70:2);
    \draw[-latex] (70:2) -- (70:3);
    \draw[thick, dotted] (70:3) -- (70:3.8);
    \draw[thick, dotted, shorten <=0pt, shorten >=0pt] (105:1.7) arc(105:85:1.7);
    \draw (-.7,-.7) node{$\mo{B}'$};
  \end{tikzpicture}
  \caption{Two structures.}
  \label{fig:two-struc}
\end{figure}

\begin{exm}
  Consider the two structures as in Figure \ref{fig:two-struc} where the lines indicate the
  relations $R$ and $R'$. Equip both models with the trivial order, that is, $x \leq y$
  iff $x = y$. Then $\mo{B}$ and $\mo{B}'$ are two strictly condensed $\Lint_{\Box}$-frames.
  
  Since the orders are taken to be trivial, the interpretation of
  intuitionistic logic is classical, i.e., every subset of states is an 
  interpretant and the interpretation of $\neg\phi$ is given by taking complements.
  Moreover, the notion of an $\Lint_{\Box}$-bisimulation reduces to a
  Kripke bisimulation in the
  usual sense for normal modal logic, see e.g.~\cite[Definition 2.16]{BRV01}.
  Therefore, the argument in Example 2.23 of {\it op.~\!cit.}~proves that
  the roots of the two models are logically equivalent but not bisimilar.
\end{exm}

\subsection{Hennessy-Milner Property for Modal Dual- and Bi-Intuitionistic Logic}

  \noindent
  We now dualise the result of Theorem \ref{thm:hm-int-box} in a similar way as
  in the proof of Theorem \ref{thm:hm-d-int} in order to obtain a Hennessy-Milner theorem
  for $\Ldint_{\Diamond}$ interpreted in $\Ldint_{\Diamond}$-models.
  This then leads to the general objective of a general Hennessy-Milner theorem for
  bi-intuitionistic modal logic with $n$ boxes and $m$ diamonds.
  
  We commence by extending Definition \ref{def:ikm-dual}, Lemma \ref{lem:sem-trans}
  and the translation $(\cdot)^t$
  to the context of \emph{modal} bi-intuitionistic logic.
  Extend the involution $(\cdot)^t$ on $\Lbiint$ to an involution on
  $\Lbiint_{\Box\Diamond}$ by adding to the recursive definition:
  $$
    (\Box\phi)^t = \Diamond\phi^t,
    \qquad
    (\Diamond\phi)^t = \Box\phi^t.
  $$
  This is easily seen to restrict to bijections
  $(\cdot)^t : \Lint_{\Box} \to \Ldint_{\Diamond}$
  and $(\cdot)^t : \Lbiint_{\Box} \to \Lbiint_{\Diamond}$.
  Furthermore, for $Z \in \{ R, S \}$ we define the \emph{dual} of a
  $Z$-model $\mo{M} = (X, \leq, Z, V)$
  to be $\mo{M}^{\partial} = (X, \geq, Z, \dual{V})$,
  where $\dual{V}(p) = X \setminus V(p)$,
  for $p \in \Prop$. Then $\mo{M}^{\partial\partial} = \mo{M}$,
  and moreover we have:
  
\begin{lem}\label{lem:1}
  The tuple $\mo{M} = (X, \leq, Z, V)$ is a (strictly condensed) $\lan{L}_{\Box}$-model
  if and only if $\dual{\mo{M}}$ is a (strictly condensed) $\lan{L}_{\Diamond}$-model.
\end{lem}

  \noindent
  Models and their duals are related in the following manner.
  This extends Lemma \ref{lem:sem-trans}.

\begin{lem}\label{lem:3}
  Let $\mo{M} = (X, \leq, R, V)$ be a strictly condensed $\lan{L}_{\Box}$-model and
  $\phi \in \Lbiint_{\Box}$ a formula. Then we have:
  $$
    \mo{M}, x \Vdash \phi \iff \mo{M}^{\partial}, x \not\Vdash \phi^t.
  $$
\end{lem}

  \noindent
  We have now set ourselves up for the proof of the Hennessy-Milner theorem
  of dual-intuitionistic logic with an extra diamond-modality.

\begin{thm}\label{thm:hm-dual-int-dia}
  Let $\mo{M} = (X, \leq, S, V)$ and $\mo{M}' = (X', \leq', S', V')$
  be two strictly condensed $\Ldint_{\Diamond}$-models such that
  $\leq$ and $\leq'$ are pre-image-compact
  and $S$ and $S'$ are image-compact.
  Then for all $x \in X$ and $x' \in X'$ we have
  $$
    x \bisim_{\Ldint_{\Diamond}} x'
    \iff x \logeq_{\Ldint_{\Diamond}} x'.
  $$
\end{thm}
\begin{proof}
  Let $B \subseteq X \times X'$ be
  the relation of logical equivalence.
  Then $B$ is also logical equivalence of $\Lint_{\Box}$-formulae
  between $\mo{M}^{\partial}$ and $(\mo{M}')^{\partial}$.
  By assumption all relations in these dual models are image-compact,
  so it follows form Theorem \ref{thm:hm-int-box} that $B$ is an
  $\Lint_{\Box}$-bisimulation between $\mo{M}^{\partial}$ and $(\mo{M}')^{\partial}$.
  An easy verification then shows that $B$ is
  an $\Ldint_{\Diamond}$-bisimulation between $\mo{M}$ and $\mo{M}'$.
\end{proof}

  \noindent
  Finally, we attain a Hennessy-Milner theorem for 
  the modal bi-intuitionistic logic $\Lbiint_{n,m}$ interpreted in $\lan{L}_{n,m}$-models.
  This follows from Theorems \ref{thm:hm-int-box} and \ref{thm:hm-dual-int-dia},
  using Lemma \ref{lem:inst} in a similar way as in the proof of
  Theorem \ref{thm:hm-bi-int}.

\begin{thm}\label{thm:hm-mod-bi-int}
  Let $\mo{M} = (X, \leq, (R_i), (S_j), V)$
  and $\mo{M}' = (X', \leq', (R_i'), (S_j'), V')$
  be two strictly condensed $\lan{L}_{n,m}$-models.
  Furthermore assume that all relations (including $\leq$ and $\leq'$) are image-compact and
  additionally that $\leq$ and $\leq'$ are pre-image-compact.
  Then for all $x \in X$ and $x' \in X'$ we have
  $$
    x \bisim_{\Lbiint_{n,m}} x'
    \iff
    x \logeq_{\Lbiint_{n,m}} x'.
  $$
\end{thm}

  \noindent
  A counterexample for the failure of the converse is readily 
  constructed from the frames $\mo{B}$ and $\mo{B}'$ in Figure
  \ref{fig:two-struc}, equipped with trivial orders $\leq$ and $\leq'$,
  and where $n = m = 1$ and $R = S$ is given by the edges, that, in general,
  logical equivalence between modal models does not imply bisimilarity.

\section{Applications}\label{sec:app}

  \noindent
  We investigate several (bi-)intuitionistic modal logics found in
  the literature, 
  and equip them with a notion of bisimulation accompanied by a Hennessy-Milner theorem.
  
  We consider (descriptive) $\Box$-models for
  the language $\Lint_{\Box}$ introduced in \cite{WolZak98} in
  Section \ref{subsec:WZ}, and in Section \ref{subsec:BD} we
  look at various ways of interpreting $\Lint_{\Box\Diamond}$ with
  a single relation for $\Box$ and $\Diamond$ (in contrast to the approach taken in
  Section \ref{sec:modal}, where each modality is interpreted via its own relation).
  In particular, this includes the well-known semantics for modal intuitionistic logic
  given by Fischer Servi \cite{Fis81},  and Plotkin and Stirling \cite{PloSti86}.

  In Subsection \ref{subsec:epistemic} we apply our results to intuitionistic epistemic
  logic \cite{JagMar16}. The knowledge operators in this logic behave like $\Box$-modalities.
  Additionally, the logic has a unary ``common knowledge'' operator $\ms{C}$,
  which behaves differently.
  
  The second half of this section is devoted to tense bi-intuitionistic logic. 
  In Subsections \ref{subsec:tense1}, \ref{subsec:tense2} and \ref{subsec:tense-H}
  we investigate three different ways of defining its semantics. 
  The corresponding notion of bisimulation requires the relations interpreting the modalities
  to look both forward and backwards. In each of these cases, we give a Hennessy-Milner class.

\subsection{Wolter/Zakarhyashev Models}\label{subsec:WZ}

  \noindent
  In \cite{WolZak98}, the authors introduce $\Box$-models as a semantics for
  $\Lint_{\Box}$. These coincide with general strictly condensed $\lan{L}_{\Box}$
  in the sense of Definition \ref{def:modal-model} with the additional property that
  the underlying order is
  a partial order (rather than a pre-order). That is:
  
\begin{defn}
  A $\Box$-frame is a tuple $(X, \leq, R, A)$ such that
  \begin{itemize}
    \item $(X, \leq)$ is a partially ordered set;
    \item $R \subseteq X \times X$ is a relation satisfying
          $({\leq} \circ R \circ {\leq}) = R$;
    \item $A \subseteq \fun{Up}(X, \leq)$ is a collection of upsets
          containing $\emptyset$ and $X$ which is closed under
          $\cap, \cup, \undto$ and $\cbox$ (cf.~Definition \ref{def:modal-general}).
  \end{itemize}
  A $\Box$-frame is called \emph{descriptive} if
  $(X, \leq, A)$ is a descriptive intuitionistic Kripke frame \cite[Section 8.4]{ChaZak97} and
  $$
    xRy \iff \forall a \in A (x \in \cbox a \text{ implies } y \in a).
  $$
  A $\Box$-model is a $\Box$-frame together with an admissible valuation
  $V : \Prop \to A$ of the proposition letters.
\end{defn}

  \noindent
  Formulae in $\Lint_{\Box}$ are interpreted as usual. 
  Since $\Box$-models are simply special cases of strictly condensed
  $\Lint_{\Box}$-models, we already have a truth-preserving notion
  of bisimulation. Moreover, Theorem \ref{thm:hm-int-box} gives rise to a
  Hennessy-Milner theorem for $\Box$-models, where image-compactness is
  now taken with respect to the general frame structure encompassed in the 
  definition of a $\Box$-model.

\begin{cor}
  Let $x$ and $x'$ be two states in two $\Box$-models all of whose
  relations are image-compact. Then
  $x \bisim_{\Lint_{\Box}} x'$ if and only if $x \logeq_{\Lint_{\Box}} x'$.
\end{cor}

  \noindent
  In particular, this holds for all descriptive $\Box$-models.
  
\begin{propn}
  Let $\mo{M} = (X, \leq, R, A)$ be a descriptive $\Box$-frame.
  Then $R$ is image-compact.
\end{propn}
\begin{proof}
  The descriptive intuitionistic Kripke frame $(X, \leq, A)$ underlying $\mo{M}$
  can be viewed as an Esakia space $(X, \leq, \tau_A)$, where $\tau_A$ is the patch
  topology defined in Definition \ref{def:general} \cite{Esa74}.
  In particular this means that $(X, \tau_A)$ is a compact topological space.
  By definition, for any $x \in X$ the set $\{ y \in X \mid x \leq y \}$
  is closed in $\tau_A$, so $\leq$ is image-compact.
  Furthermore, by definition of a descriptive $\Box$-frame we have
  $R[x] = \bigcap \{ a \in A \mid x \in \cbox a \}$
  and since this is the intersection of clopen sets, it is closed in $\tau_A$,
  hence compact.
\end{proof}

\subsection{Bo\v{z}i\'{c}/Do\v{s}en Models}\label{subsec:BD}

  \noindent
  In \cite{BozDos84}, the authors define a \emph{$\Box\Diamond$-model} to be a
  strictly condensed $\Lint_{\Box}$-model $(X, \leq, R)$ which is simultaneously an
  $\Lint_{\Diamond}$-model. That is, $(X, \leq)$ is a pre-order and $R$ is a relation
  on $X$ that satisfies
  $({\leq} \circ R \circ {\leq}) = R$ and $({\geq} \circ R) \subseteq (R \circ {\geq})$.
  These are used to interpret $\lan{Int_{\Box\Diamond}}$-formulae
  in the usual way.

  It is straightforward to see that an $\Lint_{\Box}$-bisimulation
  between $\Box\Diamond$-models preserves all formulae in $\lan{Int_{\Box\Diamond}}$,
  in particular also those involving $\Diamond$. Thus, if $x$ and $x'$ are
  two states in two $\Box\Diamond$-models with all image-compact relations,
  then we have a chain of implications:
  $$
    x \bisim_{\Lint_{\Box}} x'
    \quad\Rightarrow\quad
    x \logeq_{\Lint_{\Box\Diamond}} x'
    \quad\Rightarrow\quad
    x \logeq_{\Lint_{\Box}} x'
    \quad\Rightarrow\quad
    x \bisim_{\Lint_{\Box}} x'.
  $$
  This implies:
  
\begin{cor}
  Let $x$ and $x'$ be two states in two $\Box\Diamond$-models
  with all image-compact relations. Then
  $x \logeq_{\Lint_{\Box\Diamond}} x'$ if and only if $x \bisim_{\Lint_{\Box}} x'$.
\end{cor}

  \noindent
  We note that $\Box\Diamond$-models are special cases of the models
  used by e.g.~Fischer Servi and Plotkin and Sterling to interpret
  $\Lint_{\Box\Diamond}$, see \cite[Section 1]{PloSti86} and \cite[Section 2]{Fis81}.
  We refer to these models as FS-models, introduced formally next.
  
\begin{defn}
  An FS-model is a tuple $\mo{M} = (X, \leq, R, V)$ consisting of an intuitionistic Kripke model
  $(X, \leq, V)$ and a relation $R \subseteq X \times X$ that satisfies
  $(R \circ {\leq}) \subseteq ({\leq} \circ R)$
  and $({\geq} \circ R) \subseteq (R \circ {\geq})$.
\end{defn}
  
  \noindent
  In such a model, the interpretation of intuitionistic connectives and $\Diamond$ is
  as usual.
  However, if we interpret $\Box\phi$ as in Definition \ref{def:modal-model}
  we are no longer guaranteed an upset in $(X, \leq)$.
  This is remedied by putting
  $$
    \mo{M}, x \Vdash \Box\phi \iff \text{for all $y \in X$, } x({\leq} \circ R)y \text{ implies } \mo{M}, y \Vdash \phi.
  $$
  In the special case where
  \begin{equation}\label{eq:FS-sc}
    ({\leq} \circ R) \subseteq R,
  \end{equation}
  the interpretation of $\Box$ coincides with the one given
  in Definition \ref{def:modal-model}, i.e., without the additional quantification
  over $\leq$ in between. Moreover, if this is the case then $(X, \leq, R, V)$ is
  a strictly condensed $\Box\Diamond$-model.
  Therefore, we call an FS-model satisfying \eqref{eq:FS-sc}
  \emph{strictly condensed}. Then we have:

\begin{cor}
  Let $x$ and $x'$ be two states in two strictly condensed FS-models
  with all image-compact relations. Then
  $x \logeq_{\Lint_{\Box\Diamond}} x'$ if and only if $x \bisim_{\Lint_{\Box}} x'$.
\end{cor}

\subsection{Intuitionistic Epistemic Logic}\label{subsec:epistemic}

  \noindent
  Intuitionistic epistemic logic describes a system of the knowledge of $n$ agents
  \cite{JagMar16}.
  The logical language used for this is $\lan{EK}$, and is constructed from propositional
  variables,
  intuitionistic connectives, and additional unary modal operators
  $\ms{K}_i$ for every $i \in \{ 1, \ldots, n \}$ and $\ms{C}$.
  The intuitive meaning of $\ms{K}_i\phi$ is ``agent $i$ knows that $\phi$''
  and $\ms{C}\phi$ means that $\phi$ is common knowledge.
  This language can be interpreted in \emph{EK-models} \cite[Definitions 2 and 3]{JagMar16}.
  We give the definition of these models in a slightly reformulated way,
  so that the connection with $\Lint_{\Box}$-models is easier to see.

\begin{defn}
  An \emph{EK-model} is a tuple $(X, \leq, R_1, \ldots, R_n, V)$ consisting of
  an intuitionistic Kripke model $(X, \leq, V)$ and relations $R_i \subseteq X \times X$
  satisfying $({\leq} \circ R_i) \subseteq R_i$.
  
  The interpretation of intuitionistic connectives is as usual, and the
  interpretation of $\ms{K}_i$ is as for boxes:
  $$
    \mo{M}, x \Vdash \ms{K}_i\phi \iff \text{for all $y \in X$, } xR_iy \text{ implies } \mo{M}, y \Vdash \phi.
  $$
  The interpretation $\ms{C}$ is best described via a new relation $R^*$.
  Let $R = R_1 \cup \cdots \cup R_n$ and let $R^*$ be the collection of all pairs
  $(x, y)$ such that $y$ is reachable from $x$ via a finite number of $R$-transitions.
  Then
  $$
    \mo{M}, x \Vdash \ms{C}\phi \iff \text{for all $y \in X$, } xR^*y \text{ implies } \mo{M}, y \Vdash \phi.
  $$
\end{defn}

  Of course, EK-models are special cases of $\Lint_{n,0}$-models and the
  interpretation of the $\ms{K}_i$ corresponds to the $n$ boxes in such a model.
  A straightforward verification shows that $\Lint_{n,0}$-bisimulations
  also preserve the operator $\ms{C}$, so that we have:

\begin{lem}
  Let $x$ and $x'$ be two states in two EK-models $\mo{M}$ and $\mo{M}'$
  which are linked by an $\Lint_{n,0}$-bisimulation.
  Then $x \logeq_{\lan{EK}} x'$, that is, $x$ and $x'$ satisfy precisely the same
  $\lan{EK}$-formulae.
\end{lem}

  \noindent
  Conversely, if two states $x$ and $x'$ in two EK-models are logically equivalent,
  then in particular they satisfy the same $\Lint_{n,0}$-formulae, i.e.,
  we have $x \logeq_{\Lint_{n,0}} x'$. If $\mo{M}$ and $\mo{M}'$ (viewed as
  $\Lint_{n,0}$-models) are strictly condensed and all their relations
  are image-compact, then it follows from Theorem \ref{thm:hm-int-box} and
  Lemma \ref{lem:inst} that $x$ and $x'$ are linked by an
  $\Lint_{n,0}$-bisimulation. By the previous lemma this in turn implies
  $x \logeq_{\lan{EK}} x'$. Thus we have:

\begin{cor}
  Let $x$ and $x'$ be two states in two strictly condensed EK-models
  all of whose relations are image-compact. Then
  $$
    x \logeq_{\lan{EK}} x' \iff x \bisim_{\Lint_{n,0}} x'.
  $$
\end{cor}

  \noindent
  Therefore $\Lint_{n,0}$-bisimulations provide a suitable notion of bisimulation
  between EK-models.

\subsection{Tense Bi-Intuitionistic Logic in Tense Models}\label{subsec:tense1}\label{subsec:tense}

  \noindent
  Tense bi-intuitionistic logic is obtained from the modal bi-intuitionistic logic
  $\Lbiint_{\Box\Diamond}$ by extending it with tense operators $\tdiamond, \tbox$
  corresponding to $\Box$ and $\Diamond$, respectively. 
  We call this language $\Tbiint = \Lbiint_{2,2}$.
  Classically, $\tdiamond$ is interpreted using the converse relation
  of $\Box$. Since we assume no connection between $\Box$ and $\Diamond$, we get
  an additional tense operator $\tbox$ which is interpreted using the converse relation
  of $\Diamond$.

  In this subsection
  we adapt $\Lbiint_{\Box\Diamond}$-models (Definition \ref{def:modal-model}) to
  allow interpretation
  of $\Tbiint$-formulae, i.e., we make sure that the truth-set of every formula is still
  up-closed. In the next two subsections we investigate two more ways
  to define semantics for tense bi-intuitionistic logic.
  If $R$ is a relation on $X$, we write $\breve{R} = \{ (x,y) \mid yRx \}$
  for the converse relation.

  Let $(X, \leq, R, S, V)$ be a $\Lbiint_{\Box\Diamond}$-model
  for $\Lbiint_{\Box\Diamond}$.
  As stated, we want to use the converse relations $\breve{S}$ and $\breve{R}$ to interpret
  $\tbox$ and $\tdiamond$, respectively. Therefore, a possible semantics for
  $\Tbiint$ is given by $\Lbiint_{2,2}$-models  ${(X, \leq, R_1, R_2, S_1, S_2, V)}$
  that satisfy $R_2 = \breve{S}_1$ and $S_2 = \breve{R}_1$.
  This identification leads to the additional coherence conditions 
  $({\geq} \circ \breve{R}_1) \subseteq (\breve{R}_1 \circ {\geq})$,
  and similarly for $S_1$. Thus, we can also view such a model
  as a $\Lbiint_{1,1}$-model with additional coherence conditions.
  This is reflected in the following definition of a \emph{tense model}.

\begin{defn}\label{def:tense-model}
  A \emph{tense model} is a tuple $(X, \leq, R, S, V)$
  consisting of an intuitionistic Kripke model $(X, \leq, V)$ and
  two relations $R, S \subseteq X \times X$ satisfying
  $$
    ({\leq} \circ R) = (R \circ {\leq})
    \quad\text{and}\quad
    ({\geq} \circ S) = (S \circ {\geq}).
  $$
  \end{defn}

  \noindent
  The interpretation of the tense operators in a tense model
  $\mo{M} = (X, \leq, R, S, V)$ is given by
  \begin{align*}
    \mo{M}, x \Vdash \tdiamond\phi
      &\iff \mo{M}, y \Vdash \phi \text{ for some $y$ with } yRx \\
    \mo{M}, x \Vdash \tbox\phi
      &\iff \text{for all $y \in X$, } ySx \text{ implies } \mo{M}, y \Vdash \phi.
  \end{align*}
  Note that this corresponds precisely to the usual interpretation of box and diamond
  in the $\Lbiint_{1,1}$-model $(X, \leq, \breve{S}, \breve{R}, V)$.
  As a consequence, persistence still holds, i.e., the truth-set of every formula
  is up-closed in $(X, \leq)$.
  
  To define a \emph{tense bisimulation} between two tense models
  $(X, \leq, R, S, V)$ and $(X', \leq', R', S', V')$
  we simply use the notion of a $\Lbiint_{2,2}$-bisimulation
  between $(X, \leq, R, \breve{S}, S, \breve{R}, V)$
  and $(X', \leq', R', \breve{S}', S', \breve{R}', V')$
  from Definition \ref{def:modal-bis}.
  Explicitly, this can be defined as follows:

\begin{defn}\label{def:tense-bisim}
  By a \emph{tense bisimulation} between two tense models
  $\mo{M} = (X, \leq, R, S, V)$ and
  $\mo{M}' = (X', \leq', R', S', V')$
  we mean a $\Lbiint$-bisimulation $B \subseteq X \times X$ between the underlying
  intuitionistic Kripke models such that
  for all $(x, x') \in B$ and $Z \in \{ R, \breve{S}, S, \breve{R} \}$ we have:
  \begin{itemize}
    \item If $xZy$ then there exists $y' \in X'$ such that $x'Z'y'$ and $yBy'$;
    \item If $x'Z'y'$ then there exists $y \in X$ such that $xZy$ and $yBy'$.
  \end{itemize}
  The notion of \emph{tense bisimilarity} is defined as usual, and denoted $\bisim_{\Tbiint}$.
\end{defn}

  \noindent
  It follows from Proposition \ref{prop:adeq-modal} that $\Tbiint$-bisimilar states
  satisfy the same $\Tbiint$-formulae.

  We call a tense model $(X, \leq, R, S, V)$ \emph{strictly condensed} if 
  $({\leq} \circ R \circ {\leq}) \subseteq R$ and $({\geq} \circ S \circ {\geq}) \subseteq S$.
  A straightforward verification shows that this is the case
  if and only if
  $({\leq} \circ \breve{S} \circ {\leq}) \subseteq \breve{S}$ and
  $({\geq} \circ \breve{R} \circ {\geq}) \subseteq \breve{R}$,
  so a tense model is strictly condensed if and only if the $\Lbiint_{2,2}$-model
  $(X, \leq, R, \breve{S}, S, \breve{R}, V)$ is strictly condensed in the sense
  of Definition \ref{def:modal-model}.
  We define (pre-)image-compactness of relations in a tense model
  $(X, \leq, R, S, V)$ as if it were a $\Lbiint_{\Box\Diamond}$-model.
  As a corollary of Theorem \ref{thm:hm-mod-bi-int} we then obtain:

\begin{cor}\label{cor:hm-tense}
  Let $\mo{M}$ and $\mo{M'}$ be strictly condensed tense models all of whose 
  relations are both image-compact and pre-image-compact.
  Suppose $x \in \mo{M}$ and $x' \in \mo{M}'$. Then
  $$
    x \logeq_{\Tbiint} x' \iff x \bisim_{\Tbiint} x'.
  $$
\end{cor}

  \noindent
  We leave the construction of counterexamples showing that we cannot drop the
  conditions of (pre-)image-compactness of the relations in Corollary \ref{cor:hm-tense}
  to the reader.

\subsection{Tense Bi-Intuitionistic Logic by Gor\'{e}, Postniece and Tiu}\label{subsec:tense2}

  \noindent
  An alternative semantics for $\Tbiint$ 
  is introduced in \cite[Section 6]{GorPosTiu10}.
  The authors define a model, which we shall refer to as a
  \emph{GPT-model}, to be a tuple
  $(X, \leq, R, S, V)$ such that $(X, \leq, V)$ is an
  intuitionistic Kripke model and $R, S$ are relations
  on $X$ satisfying
  \begin{equation}\label{eq:a-tense}
    (R \circ {\leq}) \subseteq ({\leq} \circ R)
    \quad\text{and}\quad
    ({\geq} \circ S) \subseteq (S \circ {\geq}).
  \end{equation}
  The interpretation of the modalities is then given by
  \begin{align*}
    \mo{M}, x \Vdash \Box\phi
      &\iff \text{for all $y \in X$, } x({\leq} \circ R)y \text{ implies } \mo{M}, y \Vdash \phi \\
    \mo{M}, x \Vdash \Diamond\phi
      &\iff \text{there exists } y \in X
            \text{ such that } xSy
            \text{ and } \mo{M}, y \Vdash \phi \\   
    \mo{M}, x \Vdash \tbox\phi
      &\iff \text{for all $y \in X$, } x({\leq} \circ \breve{S})y \text{ implies } \mo{M}, y \Vdash \phi \\
    \mo{M}, x \Vdash \tdiamond\phi
      &\iff \text{there exists } y \in X
            \text{ such that } x\breve{R}y
            \text{ and } \mo{M}, y \Vdash \phi
  \end{align*}
  We can define a bisimulation between such models in the same way as
  in Definition \ref{def:tense-bisim} above.
  They are easily seen to preserve truth,
  despite the changed interpretation of the $\Box$-modalities.
  If a GPT-model $\mo{M} = (X, \leq, R, S, V)$
  satisfies
  \begin{equation}\label{eq:a-tense-con}
    ({\leq} \circ R) \subseteq R,
    \qquad
    ({\leq} \circ \breve{S}) \subseteq \breve{S}
  \end{equation}
  then the interpretation of $\Box$ and $\tbox$ is the same as in
  Subsection \ref{subsec:tense}, i.e., a state satisfies $\Box\phi$ (resp.~$\tbox\phi$)
  if all $R$-successors (resp.~$\breve{S}$-successors) satisfy $\phi$.
  A GPT-model that satisfies \eqref{eq:a-tense-con} will be called 
  \emph{strictly condensed}. Indeed, these are strictly condensed frames in
  the sense of Subsection \ref{subsec:tense} above, because
  \begin{align*}
    ({\leq} \circ R \circ {\leq})
      &\subseteq ({\leq} \circ {\leq} \circ R)
          &\text{(By \eqref{eq:a-tense})} \\
      &\subseteq ({\leq} \circ R)
          &\text{(${\leq}$ is transitive)} \\
      &\subseteq R
          &\text{(By \eqref{eq:a-tense-con})}
  \end{align*}
  and similarly $({\geq} \circ S \circ {\geq}) \subseteq S$.
  Since furthermore the interpretation of formulae is the same as for tense models,
  Corollary \ref{cor:hm-tense} now carries over to:

\begin{cor}
  Let $\mo{M}$ and $\mo{M'}$ be strictly condensed GPT-models all of whose 
  relations are both image-compact and pre-image-compact.
  Suppose $x \in \mo{M}$ and $x' \in \mo{M}'$. Then logical equivalence
  implies tense bisimilarity.
\end{cor}
  
  \noindent 
  As is the case for $\lan{L}_{n,m}$-models (see Proposition \ref{prop:strictify}),
  we can turn every GPT-model
  into a strictly condensed one by only modifying the relations $R$ and $S$.

\begin{propn}
  For every GPT-model $\mo{M} = (X, \leq, R, S, V)$
  we can find a strictly condensed GPT-model
  $\mo{M}^+ = (X, \leq, R^+, S^+, V)$ whose
  underlying intuitionistic Kripke model remains unchanged and which satisfies
  for all $x \in X$ and $\phi \in \Tbiint$:
  $$
    \mo{M}, x \Vdash \phi
    \iff
    \mo{M}^+, x \Vdash \phi.
  $$
\end{propn}
\begin{proof}
  Define $R^+ = ({\leq} \circ R)$ and $S^+ = (S \circ {\geq})$.
  Then reflexivity and transitivity of $\leq$ prove $({\leq} \circ R^+) = R^+$
  and $(R^+ \circ {\leq}) \subseteq ({\leq} \circ R^+)$. Besides,
  $({\geq} \circ S^+) \subseteq (S^+ \circ {\geq})$, and
  clearly $S^+ = (S \circ {\geq})$ implies
  $({\leq} \circ \breve{S}^+) \subseteq \breve{S}^+$.
  So $\mo{M}^+$ is indeed a strictly condensed GPT-model.
  
  We now prove that the theory of the individual states is unchanged,
  by induction on the structure of $\phi$. The only non-trivial cases are 
  the ones involving the modalities. We show the cases $\Box\phi$ and
  $\Diamond\phi$. Their tense counterparts are similar. We have:
  \begin{align*}
    \mo{M}, x \Vdash \Box\phi
      &\iff x({\leq} \circ R)y \text{ implies } \mo{M}, y \Vdash \phi \\
      &\iff x({\leq} \circ {\leq} \circ R)y \text{ implies } \mo{M}, y \Vdash \phi \\
      &\iff x({\leq} \circ R^+y) \text{ implies } \mo{M}^+, y \Vdash \phi \\
      &\iff \mo{M}^+, x \Vdash \Box\phi
  \end{align*}
  For the diamonds:
  \begin{align*}
    \mo{M}, x \Vdash \Diamond\phi
      &\iff \text{there exists } y \in X
            \text{ such that } xSy
            \text{ and } \mo{M}, y \Vdash \phi \\
      &\iff \text{there exists } y \in X
            \text{ such that } x(S \circ {\geq})y
            \text{ and } \mo{M}, y \Vdash \phi \\
      &\iff \text{there exists } y \in X
            \text{ such that } xS^+y
            \text{ and } \mo{M}^+, y \Vdash \phi \\
      &\iff \mo{M}^+, x \Vdash \Diamond\phi
  \end{align*}
  The second ``iff'' holds by persistence: the direction from left to right is
  immediate, conversely, if $xSz \geq y$ and $\mo{M}, y \Vdash \phi$,
  then persistence implies $\mo{M}, z \Vdash \phi$.
\end{proof}

\subsection{Tense Bi-Intuitionistic Logic in $H$-Models}\label{subsec:tense-H}

  \noindent
  Lastly, we review another approach, taken in \cite{SteSchRye16,SanSte17},
  where the authors assume additional
  axioms relating $\Box$ and $\Diamond$. In particular, in their semantics
  $\Diamond\phi$ is equivalent to $\gen\Box\neg\phi$, where $\neg\phi = \phi \to \bot$
  and $\gen\phi = \top \bito \phi$.
  The interpreting structures they use are
  \emph{$H$-frames} \cite[Definition 10]{SteSchRye16}.
  These are precisely strictly condensed $\Box$-frames from \cite{BozDos84}, called
  strictly condensed $\lan{L}_{\Box}$-frames in our notation
  (cf.~Definition \ref{def:modal-model}).
  We view an $H$-frame ($H$-model) as a strictly condensed
  $\Lbiint_{\Box}$-frame ($\Lbiint_{\Box}$-model).

  Let $\mo{M} = (X, \leq, R, V)$ be an $H$-model. 
  While $\Box$ and $\tdiamond$ are interpreted in the same way as in Subsection
  \ref{subsec:tense}, the interpretation of
  $\tbox$ and $\Diamond$ is given via the so-called \emph{left converse} of $R$,
  defined as ${\geq} \circ R \circ {\geq}$.
  Writing (suggestively)
  $S = ({\geq} \circ R \circ {\geq})$,%
  \footnote{This is the converse of $\curvearrowleftup R$ in \cite{SteSchRye16},
            which may seem odd. But verifying
            $\llb \Diamond\phi \rrb
              = \llb \phi \rrb \oplus \curvearrowleftup R
              = \{ x \in X \mid \exists y : y(\curvearrowleftup R)x
                  \text{ and } y \in \llb \phi \rrb \}
              = \{ x \in X \mid \exists y : y({\leq} \circ \breve{R} \circ {\leq})x
                  \text{ and } y \in \llb \phi \rrb \}
              = \{ x \in X \mid \exists y : x({\geq} \circ R \circ {\geq})y
                  \text{ and } y \in \llb \phi \rrb \}
            $
            shows that this is indeed how we interpret $\Diamond$.
            A similar verification shows that we get the correct interpretation for $\tbox$.}
  these modalities
  are again interpreted as usual, i.e., via
  \begin{align*}
    \mo{M}, x \Vdash \Diamond\phi &\iff \text{there exists } y \in X
                                        \text{ such that } xSy
                                        \text{ and }\mo{M}, y \Vdash \phi \\
    \mo{M}, x \Vdash \tbox\phi &\iff ySx \text{ implies } \mo{M}, y \Vdash \phi
  \end{align*}
  Therefore, setting $\ov{\mo{M}} = (X, \leq, R, S, V)$ yields
  a (strictly condensed) tense model $\ov{\mo{M}}$ in the sense of
  Definition \ref{def:tense-model} which satisfies
  $\mo{M}, x \Vdash \phi$ iff $\ov{\mo{M}}, x \Vdash \phi$.
  To see that $\ov{\mo{M}}$ is strictly condensed, note that we have
  ${\leq} \circ R \circ {\leq} = R$ by definition, and it follows from reflexivity
  and transitivity of $\leq$ that
  $$
    ({\geq} \circ S \circ {\geq})
      = ({\geq} \circ {\geq} \circ R \circ {\geq} \circ {\geq})
      = ({\geq} \circ R \circ {\geq})
      = S.
  $$
  
  \noindent
  The obvious notion of bisimulation between $H$-models is:
  
\begin{defn}
  An $H$-bisimulation between two $H$-models $(X, \leq, R, V)$ and $(X', \leq', R', V')$
  is a $\Lbiint$-bisimulation $B$ between the underlying intuitionistic Kripke models
  that additional is a $\Box$-zigzag and a $\tdiamond$-zigzag.
  (That is, both $R$ and $\breve{R}$ satisfy the zigzag conditions.)
  $H$-bisimilarity is denoted by $\bisim_H$.
\end{defn}

  \noindent
  In other words, $B$ is an $H$-bisimulation if and only if it is a
  $\Lbiint_{\Box\tdiamond}$-bi\-simu\-la\-tion between the $\Lbiint_{1,1}$-models
  $(X, \leq, R, \breve{R}, V)$ and $(X', \leq', R', \breve{R}', V')$.
  Besides, a straightforward verification shows that such an $H$-bisimulation between
  $\mo{M}$ and $\mo{M}'$ is also
  a tense bisimulation between $\ov{\mo{M}}$ and $\ov{\mo{M}}'$.
  Therefore, it preserves truth of all $\Tbiint$-formulae.

  For the converse, suppose $\mo{M} = (X, \leq, R, V)$ and
  $\mo{M}' = (X', \leq', R', V')$ are two $H$-models all of whose
  relations are image-compact and pre-image-compact. Then
  $(X, \leq, R, \breve{R}, V)$ and $(X', \leq', R', \breve{R}', V')$ are
  strictly condensed $\Lbiint_{1,1}$-models in the sense of Definition \ref{def:modal-model}.
  Moreover, they satisfy all preconditions of Theorem \ref{thm:hm-mod-bi-int}.
  If $x$ and $x'$ are two states
  in $\mo{M}$ an $\mo{M}'$ that satisfy the same $\Tbiint$-formulae, then
  in particular $x \logeq_{\Lbiint_{\Box\tdiamond}} x'$, so by Theorem \ref{thm:hm-mod-bi-int}
  there is a $\Lbiint_{\Box\tdiamond}$-bisimulation $B$ linking them.
  But by definition $B$ is precisely an $H$-bisimulation. Summarising:
  $$
    x \bisim_H x'
    \quad\To\quad x \logeq_{\Tbiint} x'
    \quad\To\quad x \logeq_{\Lbiint_{\Box\tdiamond}} x'
    \quad\To\quad x \bisim_H x'.
  $$
  Thus we have proved:
  
\begin{cor}\label{cor:hm-H-mod}
  Between any two $H$-models whose relations are image-compact and pre-image-compact,
  we have $x \logeq_{\Tbiint} x'$ if and only if $x \bisim_H x'$.
\end{cor}

\begin{rem}
  One might wonder why we did not employ the results from Subsection \ref{subsec:tense}
  in order to obtain a Hennessy-Milner result for $H$-models. This would require
  stipulating $S = ({\geq} \circ R \circ {\geq})$
  to be image-compact and pre-image-compact, on top of the preconditions of
  Corollary \ref{cor:hm-H-mod}.
  Indeed, it does necessarily follow from $\leq$ and $R$ being (pre-)image-compact.
  The current approach circumvents this.
\end{rem}

\section{Image-Compactness Versus Saturation}\label{sec:ic-vs-sat}
  
  \noindent
  We detail the relation between image-compactness and notions saturation for
  normal modal logic over a classical base, and for intuitionistic logic.
  

\subsection{Modal Saturation in Classical Modal Logic}

\noindent
We can interpret classical modal logic, that is, the language
$\Lint_\Box$, in $\Lint_{\Box}$-models where $\leq$ is equality, and
recover the classical semantics. 
  In
  particular, this implies that every subset
  is up-closed and intuitionistic negation is the same as classical negation.
  Indeed, such an $\Lint_{\Box}$-model is simply a Kripke model in the usual sense.
  We write $\lan{ML}$ for the language of classical normal modal logic.
  
  If the orders $\leq$ are trivial,
  then the definition of an $\Lint_{\Box}$-bisimulation reduces to a relation 
  that preserves truth of proposition letters and satisfies ($\Box$-zig) and ($\Box$-zag).
  In other words, it is a Kripke bisimulation for classical modal logic in the usual sense,
  see e.g.~\cite[Definition 2.16]{BRV01}.
  In this setting there is a well-known Hennessy-Milner result for the class of
  so-called \emph{m-saturated models} \cite[Proposition 2.54]{BRV01}.
  We recall the definition of m-saturation.

\begin{defn}\label{def:m-sat}
  Let $\mo{M} = (X, R, V)$ be a Kripke model and $a \subseteq X$.
  Then a set $\Sigma$ of formulae is called \emph{satisfiable} in
  $a$ if there exists a world $x \in a$ which satisfies each $\phi \in \Sigma$.
  A set $\Sigma$ is called \emph{finitely satisfiable} in $a$ if
  every finite subset of $\Sigma$ is satisfiable in $a$.
  The model $\mo{M}$ is called \emph{m-saturated}
  if for all $x \in X$ and $\Sigma \subseteq \lan{ML}$
  it satisfies:
  \begin{center}
    If $\Sigma$ is finitely satisfiable in the set of successors of $x$,\\
    then $\Sigma$ is satisfiable in the set of successors of $x$.
  \end{center}
\end{defn}

  \noindent
  Our results subsume the Hennessy-Milner result for m-saturated models
  in the following sense:
  a Kripke model $(X, R, V)$ is image-compact if and only if it
  is m-saturated.
  This result, together with the notion of image-compact relations for
  Kripke frames, also appears in \cite{BonKwi95}. 

\begin{propn}
  Let $\mo{M} = (X, R, V)$ be a Kripke model.
  Then $\mo{M}$ is image-compact if and only if it is m-saturated.
\end{propn}
\begin{proof}
  Let $x \in X$ and let $\Sigma$ be a set of formulae that is finitely satisfiable
  in the set $R[x]$ of $R$-successors of $x$.
  Suppose towards a contradiction that $\Sigma$ is not satisfiable in $R[x]$.
  Then for each $y \in R[x]$ there is a $\phi \in \Sigma$ such that
  $\mo{M}, y \not\Vdash\phi$, hence $\{ \llb \neg\phi \rrb^{\mo{M}} \mid \phi \in \Sigma \}$
  is an open cover of $R[x]$. Note that the truth set of every formula
  is clopen in $\tau_A$. By compactness of $R[x]$ we then find a finite
  subset $\Sigma' \subseteq \Sigma$ such that
  $R[x] \subseteq \bigcup_{\phi \in \Sigma'} \llb \neg\phi \rrb^{\mo{M}}$.
  But that implies that the finite set $\Sigma'$ is not satisfiable, a
  contradiction with the assumption that $\Sigma$ is finitely satisfiable.
  
  Conversely, suppose $\mo{M}$ is m-saturated.
  Let $A = \{ \llb \phi \rrb^{\mo{M}} \mid \phi \in \lan{ML} \}$. Then
  clearly $(X, R, A, V)$ is a general Kripke model. We prove that
  $R[x]$ is compact for every $x$. By the Alexander subbase theorem
  it suffices to prove that every open cover consisting of subbase elements
  has a finite subcover,
  and since $-A = A$ (because
  $X \setminus \llb \phi \rrb^{\mo{M}} = \llb \neg\phi \rrb^{\mo{M}}$
  by classicality) this subbase consists exclusively of
  truth-sets of formulae. So suppose 
  $R[x] \subseteq \bigcup_{\phi \in \Sigma} \llb \phi \rrb^{\mo{M}}$,
  for some set $\Sigma$ of formulae. Then clearly the set
  $\{ \neg\phi \mid \phi \in \Sigma \}$ is not satisfiable,
  hence (since $\mo{M}$ is m-saturated) there must be a finite $\Sigma' \subseteq \Sigma$
  such that $\{ \neg\phi \mid \phi \in \Sigma' \}$ is not satisfiable
  in $R[x]$. But that implies
  $R[x] \subseteq \bigcup_{\phi \in \Sigma'} \llb \phi \rrb^{\mo{M}}$,
  which gives the desired finite subcover.
\end{proof}

  \noindent
  In \cite{BezFonVen10} the collection of \emph{descriptive} Kripke models was
  identified as a Hennessy-Milner class.
  If $(X, R, A, V)$ is a descriptive Kripke model, then for all
  $(X, \tau_A)$ is a Stone space. Moreover $R[x]$
  is closed in $(X, \tau_A)$ for all $x \in X$, hence compact.
  Therefore, the Hennessy-Milner property for the collection of descriptive
  Kripke models also follows from our results.
  
  In \cite{Bou04}, Hennessy-Milner type results are formulated for so-called
  weak-strict languages. Such languages are interpreted in Kripke structures.
  One condition for obtaining such a result, is that the models be SW-saturated
  (Definition 3.5.1 and Lemma 3.5.8 in {\it op.~\!cit.}), which the prove to be
  equivalent to the customary notion of modal saturation in Proposition 3.5.2.

\subsection{Saturation for Intuitionistic Logic}

  \noindent
  In \cite{Pat97} several Hennessy-Milner properties for $\Lint$-bisimulations on
  intuitionistic Kripke models are given. The strongest of these uses the notion
  of \emph{local saturation}, an adaptation of m-saturation from Definition
  \ref{def:m-sat}.

\begin{defn}
  An intuitionistic Kripke model $\mo{M} = (X, \leq, V)$ is \emph{locally saturated}
  if for all $x \in X$ and disjoint sets of $\Lint$-formulae
  $\Theta_s, \Theta_r$ the following holds: 
  If for all finite subsets $\theta_s \subseteq \Theta_s$ and $\theta_r \subseteq \Theta_r$
  there are worlds $y, y' \in {\uparrow}_{\leq}x$ such that
  $\mo{M}, y \Vdash \bigwedge\theta_s$ and $y' \not\Vdash \bigvee\theta_r$,
  then there is a world $z \in {\uparrow}_{\leq}x$ which satisfies 
  every formula in $\Theta_s$ and refutes every formula in $\Theta_r$.
\end{defn}

  \noindent
  It is shown in \cite[Theorem 21]{Pat97} that logical equivalence on a locally saturated
  intuitionistic Kripke model implies $\Lint$-bisimilairty.
  We shall now show that
  an intuitionistic Kripke model is locally saturated if and only if
  it is image-compact. Therefore, Theorem \ref{thm:hm-int} is equivalent to
  {\it loc.~\!cit.}

\begin{propn}\label{prop:AP}
  An intuitionistic Kripke model $\mo{M} = (X, \leq, V)$ is locally saturated
  if and only if it is image-compact.
\end{propn}
\begin{proof}
  Suppose $\mo{M}$ is locally saturated and let $x \in X$.
  Define $A = \{ \llb \phi \rrb \mid \phi \in \Lint \}$.
  Then clearly $(X, \leq, A)$ is a general frame.
  We will show that every finite subcover of ${\uparrow}_{\leq}x = \{ y \in X \mid x \leq y \}$
  consisting of subbasic opens in $\tau_A$ has a finite subcover.
  By the Alexander subbase theorem this then proves that ${\uparrow}_{\leq}x$ is compact
  in the topological space $(X, \tau_A)$.
  Let
  \begin{equation}\label{eq:open-cover}
    \bigcup_{i \in I} \llb \phi_i \rrb^{\mo{M}} \cup
    \bigcup_{j \in J} (X \setminus \llb \psi_j \rrb^{\mo{M}})
  \end{equation}
  be an open cover of ${\uparrow}_{\leq}x$ and suppose towards a contradiction that it does
  not have a finite subcover. Then for every finite $I' \subseteq I$ and $J' \subseteq J$
  there exists $y \in {\uparrow}_{\leq}x$ such that
  $y \notin \bigcup_{i \in I'} \llb \phi_i \rrb^{\mo{M}}
       \cup \bigcup_{j \in J'} (X \setminus \llb \psi_j \rrb^{\mo{M}})
  $,
  i.e., $\mo{M}, y \Vdash \bigwedge_{j \in J'} \psi_j$ and
  $\mo{M}, y \Vdash \bigvee_{i \in I'} \phi_i$. Thus, setting
  $\Theta_s = \{ \psi_j \mid j \in J \}$ and $\Theta_s = \{ \phi_i \mid i \in I \}$,
  the precondition of weak saturatedness for ${\uparrow}_{\leq}x$ are is satisfied.
  However, there is no single $y \in {\uparrow}_{\leq}x$ which satisfies
  every $\psi_j \in \Theta_s$ and refutes every $\phi_i \in \Theta_r$, 
  because then $y$ would not be in the open cover in \eqref{eq:open-cover}.
  This contradicts the fact that $(X, \leq, V)$ is locally saturated.
  So the assumption that \eqref{eq:open-cover} has no finite subcover must be wrong,
  and we conclude that $\mo{M}$ is image-compact.
  
  Conversely, suppose $\mo{M}$ is not locally saturated.
  Then there exists $x \in X$ and collections of formulae $\Theta_s, \Theta_r$
  such that for all finite subsets $\theta_s \subseteq \Theta_s$ and
  $\theta_r \subseteq \Theta_r$
  we can find $y, y' \in {\uparrow}_{\leq}x$ such that $\mo{M}, y \Vdash \bigwedge \theta_s$
  and $\mo{M}, y' \not\Vdash \bigvee \theta_r$ while there is no $x$-successor which 
  satisfies all of $\Theta_s$ and refutes all formulae in $\Theta_r$. This means that
  $$
    \bigcup \{ \llb \phi \rrb^{\mo{M}} \mid \phi \in \Theta_r \}
    \cup \bigcup \{ \llb \psi \rrb^{\mo{M}} \mid \psi \in \Theta_s \}
  $$
  covers ${\uparrow}_{\leq}x$ but has no finite subcover. 
  Therefore $\mo{M}$ is not image-compact.
\end{proof}

\section{Conclusion and Further Research}\label{sec:conc}

  \noindent
  We have investigated the notion of \emph{image-compactness} and \emph{pre-image-com\-pact\-ness}
  for relational models that can be used to interpret classical, intuitionistic, 
  dual-intuitionistic and bi-intuitionistic (modal) logic.
  This notion allowed an efficient formulation of Hennessy-Milner theorems for
  Kripke-style bisimulations between such models.
  In classical modal logic and intuitionistic (non-modal) logic,
  our results match well-known Hennessy-Milner results \cite[Proposition 2.54]{BRV01},
  \cite[Theorem 21]{Pat97}, \cite[Corollary 3.9]{BezFonVen10},
  while for modal (dual- and bi-)intuitionistic logic we have
  described previously unknown Hennessy-Milner classes.
  In particular, the current approach generalises the results for
  (modal) bi-intuitionistic logic that were subject of the predecessor paper of the
  current paper \cite{GroPat19}.
  
  There are many interesting directions for further research.
  Firstly, we have not addressed intuitionistic logic enriched with a diamond-modality,
  i.e., $\Lint_{\Diamond}$, interpreted in $\Lint_{\Diamond}$-models.
  Inspection of the proof of Theorem \ref{thm:hm-int-box} shows that this
  no longer holds for diamonds. It would be interesting to investigate
  conditions for which $\logeq_{\Lint_{\Diamond}}$ implies $\bisim_{\Lint_{\Diamond}}$.
  Dually, this then gives rise to a Hennessy-Milner theorem for
  dual-intuitionistic logic with a box-modality.

  Second, there is the question on how to generalise this to $n$-ary 
  box- and diamond-like operators (see e.g.~\cite[Definition 1.23]{BRV01}).
  These are interpreted via $(n+1)$-ary relations, i.e.,
  $x \Vdash \Diamond(\phi_1, \ldots, \phi_n)$ if there exist
  $y_1, \ldots, y_n$ such that $(x, y_1, \ldots, y_n) \in S$ and 
  $y_i \Vdash \phi_i$ for all $i \in \{ 1, \ldots, n \}$.
  We expect that similar techniques as the ones presented in this paper will
  give rise to Hennessy-Milner properties for this generalisation of
  normal modal logic.
  
  Furthermore, in \cite{Dav09} intuitionistic logic is interpreted in topological spaces.
  These are then equipped with an additional relation that is used to interpret
  modalities $\Box$ and $\Diamond$ and their tense counterparts.
  In case the underlying topological space is an Alexandrov space, and hence
  corresponds to a pre-order, the intuitionistic connectives are interpreted as usual,
  and the modalities like in \cite{Fis81}.
  It would be interesting to see whether notions of (pre-)image-compactness can be
  extended to this setting, and how they correspond to the notion saturation given
  in \cite{Dav09}.
  
  Finally, we wonder whether the notion of image-compactness can be used or adapted to
  obtain Hennessy-Milner results for non-normal modal extensions of
  classical or (dual- or bi-)intuitionistic logic.
  In case of monotone modal logic over a classical base \cite{Han03,HanKup04}
  this has been done in \cite{Cel08}.
  It would be interesting to see how this generalises to
  monotone modal intuitionistic logic.
  Other interesting candidates for similar investigations are conditional logic 
  \cite{BalCin18,Wei19a,CiaLiu19} and 
  instantial neighbourhood logic \cite{BenEA17,BezEnqGro20}.

\paragraph{Acknowledgements}
  We would like to thank the anonymous referees for for the comprehensive comments and
  suggestions. Specifically, the references to related work helped embed our paper more
  closely into the body of existing research.

\section{Bibliography}
\bibliographystyle{elsarticle-num}
\bibliography{biblio.bib}{}

\end{document}